\documentclass[11pt]{article}

\usepackage{hyperref}
\hypersetup{colorlinks}
\usepackage[english]{babel}
\usepackage[a4paper]{geometry}
\usepackage{amsfonts}
\usepackage{amsthm}
\usepackage{amsmath}
\usepackage{graphicx}
\usepackage{subcaption}
\usepackage{amssymb}
\usepackage{enumerate}
\usepackage{xfrac}
\usepackage{array}
\usepackage{siunitx}
\usepackage{mathtools}
\usepackage{mathrsfs}
\usepackage{enumitem}
\setlist[enumerate]{topsep=0pt}
\setlist[itemize]{topsep=0pt}
\usepackage{float}
\usepackage[nodayofweek]{datetime} 
\usepackage[framemethod=tikz]{mdframed}
\usepackage{tikz}
\graphicspath{ {images/} }
\usepackage{listings}
\usepackage{color}
\usepackage[draft]{todonotes} 
\usepackage{comment}
\usepackage{stmaryrd}
\usepackage{bm}
\usepackage{pgfplots}
\pgfplotsset{compat=1.18}

\usepackage{nicefrac}
\usepackage{placeins}
\definecolor{dkgreen}{rgb}{0,0.6,0}
\definecolor{gray}{rgb}{0.5,0.5,0.5}
\definecolor{mauve}{rgb}{0.58,0,0.82}

\allowdisplaybreaks
\numberwithin{equation}{section} 

\newtheorem{theorem}{Theorem}[section]

\newtheorem{lemma}[theorem]{Lemma}
\newtheorem{corollary}[theorem]{Corollary}

\theoremstyle{definition}
\newtheorem{definition}[theorem]{Definition}
\newtheorem{assumption}[theorem]{Assumption}
\newtheorem{remark}[theorem]{Remark}

\newcommand{\half}{\nicefrac{1}{2}}
\newcommand{\R}{\mathbb{R}}

\newcommand{\N}{\mathbb{N}}

\newcommand{\prob}{\mathbb{P}}

\newcommand{\cP}{\mathcal{P}}
\newcommand{\cT}{\mathcal{T}}

\newcommand{\polOrder}{\ell}
\newcommand{\vemSpace}{V_{h,\polOrder}}
\newcommand{\enlargedVemSpace}{\widetilde{V}_{h,\polOrder}}

\newcommand{\element}{E}
\newcommand{\valueProj}{\Pi^{\element}_0}
\newcommand{\valueGlobal}{\Pi^h_0}
\newcommand{\edgeProj}{\Pi^e_{0}}
\newcommand{\gradProj}{\Pi^{\element}_{1}}

\newcommand{\localEnlargedVemSpace}{\enlargedVemSpace^{\element}}
\newcommand{\vertiii}[1]{{\left\vert\kern-0.25ex\left\vert\kern-0.25ex\left\vert #1 
    \right\vert\kern-0.25ex\right\vert\kern-0.25ex\right\vert}}
\newcommand{\mesh}{\cT_h}
\newcommand{\localVEM}{V_{h,\polOrder}^{\element}}

\newcommand{\dx}{\mathrm{d}\boldsymbol{x}}
\newcommand{\ds}{\mathrm{d}s}

\newcommand{\red}[1]{\textcolor{red}{#1}}

\newcommand{\edgeSet}{\mathcal{E}_h}
\newcommand{\error}{\xi}

\newcommand\normal[1][e]{\boldsymbol{n}_{#1}}

\begin{document}

\title{A posteriori error analysis of the virtual element method for second-order quasilinear elliptic PDEs}
\author{Scott Congreve and Alice Hodson\thanks{Corresponding author (\href{mailto:hodson@karlin.mff.cuni.cz}{hodson@karlin.mff.cuni.cz}).
\\
Charles University, Faculty of Mathematics and Physics, Sokolovsk\'a 83, 186 75, Praha, Czech Republic} }

\date{}

\maketitle

\begin{abstract}
    In this paper we develop a $C^0$-conforming virtual element method (VEM) for a class of second-order quasilinear elliptic PDEs in two dimensions.
We present a posteriori error analysis for this problem and derive a residual based error estimator.
The estimator is fully computable and we prove upper and lower bounds of the error estimator which are explicit in the local mesh size.
We use the estimator to drive an adaptive mesh refinement algorithm.
A handful of numerical test problems are carried out to study the performance of the proposed error indicator.
\end{abstract}


\normalsize
\section{Introduction}\label{sec: intro}

In this paper we present a conforming virtual element method (VEM) of arbitrary order for the numerical solution of a quasilinear elliptic problem in two dimensions.
We consider the \emph{a posteriori} error analysis in the $H^1$-seminorm of the virtual element method for the following boundary value problem:
\begin{equation}\label{eqn: pde}
    \begin{split}
        -\nabla \cdot \left( \mu(\boldsymbol{x},|\nabla u|)\nabla u \right) &= f \quad \text{in } \Omega, \\
        u &= 0 \quad \text{on } \partial \Omega
    \end{split}
\end{equation}
for a polygonal domain $\Omega \subset \R^2$ and $f \in L^2(\Omega)$.
For ease of presentation, we only consider homogeneous Dirichlet boundary value problems.
Furthermore, we assume that the nonlinearity $\mu$ satisfies the following set of standard assumptions see e.g. \cite{houston2008posteriori}.
\begin{assumption}[Nonlinearity assumptions]\label{ass: mu}
    We assume that the nonlinearity $\mu$ satisfies the following conditions.
    \begin{enumerate}[label=$(\alph*)$]
        \item $\mu \in C^0(\overline{\Omega} \times [0,\infty))$
        \item \label{assumption: A2} There exist positive constants $m_{\mu}, M_{\mu}$ such that
        $$ m_{\mu} (t-s) \leq \mu(\boldsymbol{x},t)t - \mu(\boldsymbol{x},s)s \leq M_{\mu} (t-s), \quad t \geq s \geq 0, \text{ for all } \boldsymbol{x} \in \overline{\Omega}.$$
    \end{enumerate}
\end{assumption}
Importantly, if $\mu$ satisfies \ref{assumption: A2} in Assumption~\ref{ass: mu}, then it can be shown that there exist constants $C_1,C_2$ with $C_1 \geq C_2 >0$ such that for any $v,w \in \R^2$ and $\boldsymbol{x} \in \overline{\Omega}$,
\begin{align}
    |\mu(\boldsymbol{x},|v|)v - \mu(\boldsymbol{x},|w|)w | &\leq C_1 |v-w|,
    \label{eqn: mu prop 1}
    \intertext{and}
    C_2 |v-w|^2 &\leq (\mu(\boldsymbol{x},|v|)v - \mu(\boldsymbol{x},|w|)w) \cdot (v-w); \label{eqn: mu prop 2}
\end{align}
cf., \cite[Lemma 2.1]{liu1994quasi}. We note that the above assumptions are fulfilled by several physical models from continuum mechanics; e.g., the Carreau law $$\mu(\boldsymbol{x},t)=k_\infty + (k_0-k_\infty)(1+\lambda t^2)^{\nicefrac{(\theta-2)}2},$$ with $k_0 > k_\infty > 0$ and $\theta\in (1,2]$.
Moreover, it is worth mentioning that studying problem \eqref{eqn: pde} is an important stepping stone in deriving efficient and effective numerical methods for non-Newtonian flow problems, see e.g. \cite{congreve2013discontinuous} and the references therein.

Introduced in 2013 \cite{beirao_da_veiga_basic_2013}, the virtual element method is an extended and generalised version of the finite element and mimetic finite difference methods and was first introduced in the context of second-order elliptic problems.
The method is highly advantageous for many reasons including the ease with which the method extends to general polygonal and polyhedral meshes.
This is extremely beneficial especially when developing adaptive schemes, since hanging nodes are automatically permissible within the VEM framework.
Furthermore, their versatility has been showcased by the wide range of problems they have been applied to over the past 10 years.
These include, but are not limited to, the following: the development of higher order continuity spaces \cite{antonietti2021review,antonietti_conforming_2018}, $H^m$-conforming VEM in any dimension \cite{chen2022conforming},
as well as the application to nonlinear problems such as the steady  Navier-Stokes equation \cite{da2018virtual} where the discrete velocity field is shown to be pointwise divergence-free.
We note that the first adaptive VEM schemes appeared in \cite{beirao2015residual,berrone2017residual,cangiani2017posteriori} and additional further developments of VEM and its applications can be found in the recent book \cite{antonietti2022virtual}.


A virtual element discretisation of a quasilinear diffusion problem in both two and three dimensions is considered in \cite{cangiani2020virtual}, while a two grid virtual element algorithm is developed in \cite{chen2024two} and a priori estimates in the $H^1$-norm are derived.
Other polygonal methods which have been considered for the discretisation of problem \eqref{eqn: pde} can be found in the following works: the mimetic finite difference (MFD) method is analysed in \cite{antonietti2015mimetic} where only a lowest order approximation is considered and no a posteriori analysis is carried out,
an HHO method for a general class of Leray-Lions elliptic equations (including the problem considered in this paper) is presented in \cite{di2017hybrid} as well as a class of quasilinear elliptic problems of nonmonotone type are analysed in \cite{gudi2022hybrid}.
Similarly to the MFD case, only a priori estimates are shown.

In this paper we present a virtual element method for problem \eqref{eqn: pde}.
We introduce a $C^0$-conforming VEM based on those introduced in \cite{ahmad_equivalent_2013,cangiani_conforming_2015} where we follow the projection approach detailed in \cite{10.1093/imanum/drab003,dedner2022framework}.
In this approach, the projection operators are defined without using the underlying variational problem, allowing us to apply the method directly to nonlinear problems including quasilinear problems such as \eqref{eqn: pde}.
We note that the same projection method has also been applied to the nonlinear fourth-order Cahn-Hilliard equation in \cite{dedner2021higher};
however, only a priori error analysis is carried out in this case.
The aforementioned projection method involves defining a hierarchy of projections, starting with a constrained least squares (CLS) problem for the value projection.
All projections are fully computable from the degrees of freedom (dofs) and are shown to be $L^2$ projections.
In this work, we employ this hierarchical projection approach in our discrete construction.
Consequently, we are able to discretise the nonlinearity $\mu$ directly using the gradient projection, which itself is an $L^2$ projection of the gradient \cite{10.1093/imanum/drab003}.
This approach avoids complicated averaging techniques seen in e.g. \cite{antonietti_$c^1$_2016} where special treatment of the nonlinearity is required and thus restricted their approach to a lowest order approximation.

Virtual element methods for quasilinear problems have been studied in \cite{cangiani2020virtual}; however, in contrast to the results shown in \cite{cangiani2020virtual}, we present an a posteriori analysis which follows the ideas introduced in \cite{cangiani2017posteriori}.
We carry out the analysis under the same regularity assumptions required in the linear setting \cite{beirao_da_veiga_basic_2013} thus allowing very general polygonal meshes, which we exploit in our adaptive algorithm.
To the best of our knowledge, this work is the first a posteriori error analysis for a virtual element discretisation of elliptic quasilinear problems on general polygonal meshes.

The structure of this paper is as follows.
In section~\ref{sec: virtual element method} we introduce the proposed virtual element method and setup the discrete projections, spaces, and forms.
Furthermore, by employing results from the theory of monotone operators, we show that the discrete problem has a unique solution.
We carry out the a posteriori error analysis in section~\ref{sec: posteriori anal}, deriving a reliable and efficient residual based error estimator.
We carry out a series of numerical experiments in section~\ref{sec: numerics} and use the estimator to drive an adaptive algorithm.
Finally, concluding remarks are given in section~\ref{sec: conclusion}.

Throughout this paper we adopt the standard notation for Sobolev spaces $H^s(\mathcal{D})$ for nonnegative integers $s$, and domains $\mathcal{D}$, with the norm and seminorm denoted by $\| \cdot \|_{s,\mathcal{D}}$ and $| \cdot |_{s,\mathcal{D}}$, respectively.
When $ \mathcal{D} = \Omega$ we may omit the subscript.

We first write problem \eqref{eqn: pde} in variational form: find $u \in H^1_0(\Omega)$ such that
\begin{align}\label{eqn: var form pde}
    a(u;u,v) = (f,v), \quad \forall v \in H^1_0(\Omega),
\end{align}
where the form $a$ is given by $a(u;v,w) := (\mu(\boldsymbol{x},|\nabla u |) \nabla v, \nabla w)$ and $(\cdot,\cdot)$ denotes the standard $L^2$ inner product over $\Omega$.
We note that, assuming sufficient regularity for the right hand side $f$, problem \eqref{eqn: pde} admits a unique solution, see e.g. \cite{douglas1971uniqueness}.

\section{The virtual element method}\label{sec: virtual element method}

In this section, we discuss the numerical approximation of problem \eqref{eqn: pde} by the virtual element method.
Following the discretisation approach in \cite{ahmad_equivalent_2013,cangiani_conforming_2015}, we use the so called ``VEM enhancement'' technique (which has been extended to fourth-order problems in \cite{10.1093/imanum/drab003}) to discretise problem \eqref{eqn: var form pde}.

\subsection{Mesh regularity and a polynomial approximation result}
Let $\mesh$ denote a tessellation of the computational domain $\Omega \subset \R^2$ into simple nonoverlapping polygons $\element$ such that $\overline\Omega = \bigcup_{\element\in\mesh} \element$. We denote by $h_{\element} := \text{diam}(\element)$ the diameter of $\element\in\mesh$ and let $h=\max_{\element\in\mesh} h_{\element}$ denote the maximum diameter of all elements. We define the set of all edges in the mesh by $\edgeSet$, which we split into boundary edges $\edgeSet^{bdy}$ and interior edges $\edgeSet^{int}$, such that
$\edgeSet^{bdy} := \{ e \in \edgeSet : e \subset \partial \Omega\}$ and $\edgeSet^{int} = \edgeSet \backslash \edgeSet^{bdy}$.

Let $\boldsymbol{v}$ be a vector-valued function, which is smooth inside each element $\element\in\mesh$.
Given two adjacent elements ${\element}^+,\element^-\in\mesh$ sharing a common edge $e\in\edgeSet^{int}$, i.e., $e \subset \partial {\element}^+ \cap \partial \element^-$, we write $\boldsymbol{v}^{\pm}$ to denote the trace of $\boldsymbol{v}|_{\element^{\pm}}$ on the edge $e$ taken from the interior of $\element^{\pm}$, respectively.
We then define the jump operator across the edge $e \in \edgeSet$ as follows: $ \llbracket \boldsymbol{v} \rrbracket := \boldsymbol{v}^{+} \cdot \normal^+ + \boldsymbol{v}^- \cdot \normal^- $ where $\normal^{\pm}$ denotes the unit outward normal on $e$ from $\element^{\pm}$, respectively.
On a boundary edge, $e \in \edgeSet^{bdy}$, we let $\llbracket \boldsymbol{v} \rrbracket := \boldsymbol{v} \cdot \normal$, where $\normal$ is the unit outward normal on $\partial\Omega$.

\begin{assumption}[Mesh regularity assumptions]\label{ass: mesh regularity}
    We assume there exists a constant $\rho >0$ such that
    \begin{enumerate}[label={$(\alph*)$}]
        \item each element of the mesh $\element \in \mesh$ is star-shaped with respect to a ball of radius $\rho h_{\element}$, and
        \item $h_e \geq \rho h_{\element}$ for every element $\element \in \mesh$ and every edge $e \subset \partial \element$.
        \label{ass: mesh reg part b}
    \end{enumerate}
    Note that these assumptions are standard in the virtual element setting as found in \cite{beirao_da_veiga_basic_2013}.
\end{assumption}
\begin{remark}\label{rmk: subtriangulation}
    Importantly, we note that as a consequence of the above, each element $\element \in \mesh$ admits a sub-triangulation; i.e., a partition of $\element$ into triangles.
    This can be obtained by joining each vertex of $\element$ to a point $\boldsymbol{x}_{\element}$ in $\element$ such that $\element$ is star shaped with respect to $\boldsymbol{x}_{\element}$.
\end{remark}
For any  $k \in \N$, we denote by $\prob_{k}(D)$ the space of polynomials of degree at most $k$ on $D\subset\R^d$, $d=1,2$. In practice, $D$ is either an element $\element\in\mesh$ or an edge $e\in\edgeSet$.
For an element $\element \in \mesh$, we denote by $\cP^{\element}_{k} : L^2(\element) \rightarrow \prob_{k}(\element)$ the $L^2$-orthogonal projection onto $\prob_{k}(\element)$.
The following theorem is an important result for the theory in section~\ref{sec: posteriori anal}, the proof of which can be obtained following the theory in either \cite{brenner_mathematical_2008,clarlet1987finite}.
\begin{theorem}[Approximation using polynomials]\label{thm: polynomial approximation}
    Under Assumption~\ref{ass: mesh regularity}, for any $k \geq 0$ and for any $w \in H^m(\element)$ with $1 \leq  m \leq k+1$, it holds that
    \begin{align*}
        \| w - \cP^{\element}_{k} w \|_{0,\element} + h_{\element} | w - \cP^{\element}_{k} w |_{1,\element} \leq C_{3} h_{\element}^m |w|_{m,\element}
    \end{align*}
    where the constant $C_{3}$ depends only on $k$ and the mesh regularity.
\end{theorem}

Furthermore, we note that we can split the form $a(\cdot;\cdot,\cdot)$ and norm $| \cdot |_{1}$ as
\begin{align*}
    a(u;v,w) := \sum_{\element \in \mesh} a^{\element}(u;v,w) &= \sum_{\element \in \mesh} (\mu(\boldsymbol{x},|\nabla u |) \nabla v, \nabla w)_{\element}, \quad &&\forall u,v,w \in H^1_0(\Omega),
    \\
    |v|_{1} &:= \left( \sum_{\element \in \mesh}  |v|_{1,\element}^2 \right)^{\half} \ &&\forall v \in H^1_0(\Omega).
\end{align*}

\subsection{The discrete spaces and projection operators}

In this section, our aim is to build a discrete VEM space $\vemSpace \subset H^1_0(\Omega)$ of polynomial order $\polOrder \geq 1$, our order of approximation, such that, for each $\element \in \mesh$, $\prob_{\polOrder}(\element) \subset \vemSpace|_{\element}$.
The other crucial aspect of the virtual element discretisation is the construction of suitable fully computable projection operators $\valueProj, \gradProj$.

To this end, we begin by building an enlarged VEM space $\localEnlargedVemSpace$ and an appropriate set of extended \emph{degrees of freedom} (dofs) for this space.
We then introduce a reduced set of dofs for $\localVEM$ and use these dofs to construct \emph{dof compatible} projection operators, which will be used to define the local VEM space as well as the discrete forms.
Further details for topics in this section can be found in \cite{10.1093/imanum/drab003}.
\begin{definition}
    For an element $\element\in\mesh$ we define the \emph{enlarged VEM space} as
    \begin{align*}
        \localEnlargedVemSpace := \left\{ v_h \in H^1(\element) : \Delta v_h \in \prob_{\polOrder}(\element) \text{ and } v_h|_{e} \in \prob_{\polOrder}(e) \, \forall e \subset \partial \element \right\},
    \end{align*}
    which we characterise by the set of extended degrees of freedom $\widetilde{\Lambda}^{\element}$ described using the \emph{dof tuple} notation introduced in \cite{10.1093/imanum/drab003}; i.e., $\widetilde{\Lambda}^{\element}$ is characterised by the dof tuple $$(0,-1,\polOrder-2,-1,\polOrder).$$
    That is, the dofs $\widetilde{\Lambda}^{\element}$ characterising the enlarged VEM space for ${v_h \in H^1(\element)}$ are:
    \begin{enumerate}[label={$(D\arabic*)$}]
        \item The value of $v_h$ at each vertex of $\element$.
        \item For $\polOrder > 1$, the moments of $v_h$ up to order $\polOrder-2$ on each edge $e \subset \partial \element$
        \begin{align*}
            \frac{1}{|e|} \int_e v_h p \, \ds \quad \forall p \in \prob_{\polOrder-2}(e).
        \end{align*}
        \item For $\polOrder > 1$, the moments of $v_h$ up to order $\polOrder$ inside $\element$
        \begin{align*}
            \frac{1}{|\element|} \int_{\element} v_h p \, \ds \quad \forall p \in \prob_{\polOrder}(\element).
        \end{align*}
    \end{enumerate}
    A proof of the unisolvency of these degrees of freedom can be found in \cite{cangiani_conforming_2015}.
\end{definition}

We can now define the local VEM space $\localVEM$, for an element $\element\in\mesh$,  as a subspace of the enlarged space $\localEnlargedVemSpace$.
In order to do so, we first introduce an interior value projection ${\valueProj : \localEnlargedVemSpace \rightarrow \prob_{\polOrder}(\element)}$ and an edge value projection $\edgeProj : \localEnlargedVemSpace \rightarrow \prob_{\polOrder}(e)$.
These projections must be computable from the reduced set of degrees of freedom $\Lambda^{\element}$, described by the dof tuple,
\begin{align}\label{eqn: reduced dof tuple}
    (0,-1,\polOrder-2,-1,\polOrder-2);
\end{align}
i.e., $\Lambda^{\element}$ consists of the original dofs introduced in e.g. \cite{beirao_da_veiga_basic_2013} for second-order problems.
We also require that these projections satisfy the following assumptions.

\begin{assumption}\label{ass: projection properties}
    For $v_h \in \localEnlargedVemSpace$, we assume that the value and edge projection operators $\valueProj$ and $\edgeProj$ are a linear combination of the degrees of freedom $\Lambda^{\element}(v_h)$.
    Furthermore, we assume they satisfy the following properties.
    \begin{enumerate}[label={$(\alph*)$}]
        \item The value projection $\valueProj v_h \in \prob_{\polOrder}(\element)$ satisfies
        \begin{align*}
            \int_{\element} \valueProj v_h p \, \dx = \int_{\element} v_h p \, \dx \quad \forall p \in \prob_{\polOrder-2}(\element)
        \end{align*}
        and $\valueProj q = q$ for all $q \in \prob_{\polOrder}(\element)$.
        \item For each edge $e \subset \partial \element$, the edge projection $\edgeProj \in \prob_{\polOrder}(e)$ satisfies $\edgeProj v_h (e^{\pm}) = v_h (e^{\pm})$, where $e^{\pm}$ denotes the vertices of an edge $e$,
        \begin{align*}
            \int_e \edgeProj v_h p \, \ds = \int_e v_h p \, \ds \quad \forall p \in \prob_{\polOrder-2}(e),
        \end{align*}
        and $\edgeProj q = q |_e$ for all $q \in \prob_{\polOrder}(\element)$.
    \end{enumerate}
\end{assumption}

We note that there are multiple ways of defining the value and edge projections such that Assumption~\ref{ass: projection properties} is satisfied.
An example choice based on a constrained least squares problem can be found in \cite{10.1093/imanum/drab003,dedner2022framework}, where the reader can find more details.
We use the choice from \cite{dedner2022framework} for the numerical experiments in section~\ref{sec: numerics}.
Now, assuming we have a value and edge projection satisfying Assumption~\ref{ass: projection properties}, we define the gradient projection and the local virtual element space $\localVEM$.

\begin{definition}
    The \emph{gradient projection} $\gradProj : \localEnlargedVemSpace \rightarrow [\prob_{\polOrder-1}(\element)]^2$ is defined as
    \begin{align*}
        \int_{\element} \gradProj v_h \cdot \boldsymbol{p} \, \dx = -\int_{\element} \valueProj v_h \nabla \cdot \boldsymbol{p} \, \dx + \sum_{e \subset \partial \element} \int_e \edgeProj v_h \boldsymbol{p} \cdot \normal \, \ds
        \quad \forall \boldsymbol{p} \in [\prob_{\polOrder-1}(\element)]^2,
    \end{align*}
    where $\normal$ denotes the unit outward normal vector to the edge $e$.
\end{definition}

\begin{definition}
    We define the \emph{local virtual element space} $\localVEM$ as
    \begin{align*}
        \localVEM := \left\{ v_h \in \localEnlargedVemSpace \, : \, (v_h - \valueProj v_h , p)|_{\element} = 0 \quad \forall p \in \prob_{\polOrder}(\element) \backslash \prob_{\polOrder-2}(\element) \right\},
    \end{align*}
    for each element $\element\in\mesh$.
\end{definition}
The local space is characterised by the dof set $\Lambda^{\element}$ described by the dof tuple in \eqref{eqn: reduced dof tuple}, and proof that this dof set is unisolvent can be found in \cite{cangiani_conforming_2015}.
Importantly, as shown in \cite{10.1093/imanum/drab003}, we have the following crucial property of our projection operators.
\begin{lemma}\label{lemma: l2 property}
    Assume that the value and edge projections satisfy Assumption~\ref{ass: projection properties}. Then, for any $v_h \in \localVEM$, it holds that
    \begin{align}
        \valueProj v_h &= \cP^{\element}_{\polOrder} v_h
        \label{eqn: L2 property value proj}
        \intertext{and}
        \gradProj v_h &= \cP^{\element}_{\polOrder-1} (\nabla v_h)
        \label{eqn: L2 property grad proj}
    \end{align}
    i.e., the value projection is the $L^2$-orthogonal projection of order $\polOrder$ and the gradient projection is the $L^2$-orthogonal projection of order $\polOrder-1$ of the gradient.
\end{lemma}

\subsection{Global spaces and the discrete forms}
With the definition of the local spaces for each element $\element\in\mesh$ in place, we can now define the global VEM space $\vemSpace$ as
\begin{align*}
    \vemSpace := \left\{ v_h \in H^1_0(\Omega) \, : \, v_h |_{\element} \in \localVEM \quad \forall \element \in \mesh \right\} \subset H^1_0(\Omega),
\end{align*}
where the global degrees of freedom are defined from the local degrees of freedom in the usual way, cf. \cite{beirao_da_veiga_basic_2013}, with local degrees of freedom corresponding to boundary vertices and boundary edges set to zero.

Lastly, we define the discrete forms required for the VEM formulation. We construct our virtual form elementwise as
\begin{align}\label{eqn: discrete form}
    a_h(z_h;v_h,w_h) = \sum_{\element \in \mesh} a_h^{\element}(z_h;v_h,w_h)
\end{align}
for any $z_h,v_h,w_h \in \vemSpace$, where
\begin{align*}
    a_h^{\element} (z_h;v_h,w_h) := (\mu(\boldsymbol{x},|\gradProj z_h |) \gradProj v_h,\gradProj w_h )_{\element} + S^{\element} (z_h ; (I-\valueProj) v_h, (I-\valueProj) w_h )
\end{align*}
for some \emph{admissible} stabilising form $S^{\element}(\cdot;\cdot,\cdot)$.
We call the stabilisation $S^{\element}(\cdot;\cdot,\cdot)$ \emph{admissible} if
    there exist positive constants $C_*,C^*$, independent of $h, \element$, such that, for all $z_h,v_h \in \localVEM$ and all $E\in\mesh$,
\begin{align}\label{eqn: admissible stabilisation}
    C_* a^{\element}(z_h;v_h,v_h) \leq S^{\element} (z_h;v_h,v_h) \leq C^* a^{\element}(z_h;v_h,v_h).
\end{align}
\begin{definition}[Stabilisation]\label{def: stabilisation}
    We define the stabilisation as
    \begin{align*}
        S^{\element} ( z_h; v_h, w_h) := M_{\mu} m_{\mu} \sum_{\lambda \in \Lambda^{\element}} \lambda(v_h)\lambda(w_h),
    \end{align*}
    making use of the standard \emph{dofi-dofi} stabilisation; cf. \cite{beirao_da_veiga_basic_2013}.
\end{definition}
Following the usual scaling argument from e.g. \cite{beirao_da_veiga_basic_2013}, it is clear that there exist positive constants $\beta_*,\beta^*$ such that for any $z_h \in \vemSpace$,
\begin{align*}
    \beta_* \int_{\element}\nabla v_h \cdot \nabla v_h \, \dx \leq S^{\element} ( z_h; v_h, v_h)  \leq \beta^* \int_{\element} \nabla v_h \cdot \nabla v_h \, \dx.
\end{align*}
From Assumption~\ref*{ass: mu}\ref{assumption: A2} with any $t>0$ and $s=0$, we have that 
\begin{align}\label{eqn: mu bound}
    m_{\mu} \leq \mu(\boldsymbol{x},t) \leq M_{\mu}.
\end{align}
Therefore, taking $t=|\nabla z_h|$ in \eqref{eqn: mu bound}, it is clear to see that \eqref{eqn: admissible stabilisation} holds with $C_* := \beta_* (M_{\mu})^{-1}$ and $C^*:= \beta^* (m_{\mu})^{-1} $. 

\begin{remark}
    Alternatively, we could follow the stabilisation approach taken in \cite{adak2022implementation,cangiani2020virtual} and define the stabilisation as
    \begin{align*}
        S^{\element} ( z_h; v_h, w_h) := \mu_{\element}(\boldsymbol{x},| \Pi^{\element,0}_{1} z_h |) \sum_{\lambda \in \Lambda^{\element}} \lambda(v_h)\lambda(w_h),
    \end{align*}
    where $\Pi^{\element,0}_{1}$ denotes the gradient projection onto the space of constant polynomials, that is, ${\Pi^{\element,0}_{1} : \localEnlargedVemSpace \rightarrow [\prob_{0}(\element)]^2}$, and $\mu_{\element}(\cdot)$ denotes the average of the function $\mu$ on the element $\element$.
    It is straightforward to show that this choice of stabilisation also satisfies \eqref{eqn: admissible stabilisation}.
    We note that this is the stabilisation we use in section~\ref{sec: numerics}; however both choices of stabilisation exhibit very similar numerical results.
    The only reason we use the linear stabilisation in Definition~\ref{def: stabilisation} throughout sections~\ref{sec: virtual element method} and \ref{sec: posteriori anal} is to ensure the validity of Lemma~\ref{lemma: lipschitz and monotone}, without requiring further restrictions on the nonlinearity $\mu$. 
\end{remark}

We also note the following important property holds. For any $v_h,z_h \in \vemSpace$
\begin{align}\label{eqn: stab bound}
    \| \nabla v_h \|_{0,\element} \leq C_4 (S^{\element}(z_h;v_h,v_h))^{\half}
\end{align}
where $C_4 :=  \min{((C_* M_{\mu})^{-\half},(C_* m_{\mu})^{-\half})}$.
Additionally, for every $\element \in \mesh$, we have the crucial stability property for any admissible stabilising form $S^{\element}$:
    there exist constants $\alpha_*, \alpha^*$, independent of $h$ and $\element$, such that for all $v_h, z_h \in \localVEM$,
    \begin{align}\label{eqn: stability}
        \alpha_* a^{\element}(z_h;v_h,v_h) \leq a_h^{\element} (z_h; v_h,v_h) \leq \alpha^* a^{\element}(z_h;v_h,v_h).
    \end{align}

\subsection{The discrete problem}
In this section we state the virtual element method for \eqref{eqn: pde} and show, in  Theorem~\ref{thm: existence and uniqueness}, that it has a unique solution.

For order $\polOrder \geq 1$, the \emph{virtual element method} discretisation of problem \eqref{eqn: var form pde} reads as follows: find $u_h \in \vemSpace$ such that
    \begin{align}\label{eqn: vem problem}
        a_h(u_h;u_h,v_h) =  L_h(v_h), \quad \forall v_h \in \vemSpace
    \end{align}
    where the right hand side $L_h(v_h) = (f_h,v_h) := \sum_{\element \in \mesh} (\valueProj f ,v_h)_{\element}$, and the discrete form $a_h(\cdot;\cdot,\cdot)$ is defined in \eqref{eqn: discrete form}.

\begin{theorem}[Existence and uniqueness of a discrete solution]\label{thm: existence and uniqueness}
    Under Assumption~\ref{ass: mu}, for a given $f \in L^2(\Omega)$ there exists a unique element $u_h \in \vemSpace$ such that \eqref{eqn: vem problem} holds.
\end{theorem}

The main tool required to prove Theorem~\ref{thm: existence and uniqueness}, is the following Lemma which details two important properties of the discrete form $a_h$,
the proof of which is a straightforward consequence of the stability properties \eqref{eqn: admissible stabilisation}, \eqref{eqn: stability} as well as properties of the nonlinearity \eqref{eqn: mu prop 1}, \eqref{eqn: mu prop 2}.

\begin{lemma}\label{lemma: lipschitz and monotone}
    The discrete form $a_h$ defined in \eqref{eqn: discrete form} admits the following properties.
    \begin{enumerate}
        \item[$(a)$]
        $a_h$ is \emph{Lipschitz continuous}, in the sense that
        \begin{align}\label{eqn: lipschitz cts}
            \left| a_h(w_h;w_h,v_h) - a_h(z_h;z_h,v_h) \right| \leq C | w_h - z_h |_{1} | v_h |_{1}
        \end{align}
        for all $w_h,z_h,v_h \in \vemSpace$.

        \item[(b)]
        $a_h$ is \emph{strongly monotone} in the sense that
        \begin{align}\label{eqn: strongly monotone}
            a_h(w_h;w_h,w_h-z_h) - a_h(z_h;z_h,w_h - z_h) \geq C | w_h - z_h |^2_{1}
        \end{align}
        for all $w_h,z_h \in \vemSpace$.
    \end{enumerate}
\end{lemma}

\color{black}
We omit the proof of Theorem~\ref{thm: existence and uniqueness} since we can directly apply results from the theory of monotone operators as detailed in \cite[Theorem 2.5]{houston2005discontinuous} with \cite[Lemma 2.2 and Lemma 2.3]{houston2005discontinuous} replaced by Lemma~\ref{lemma: lipschitz and monotone}.

\section{A posteriori error analysis}\label{sec: posteriori anal}

In this section we carry out the a posteriori error analysis for the standard $H^1$ seminorm $| \cdot|_1$; we note that due to the boundary conditions as well as the Poincar\'{e}  inequality, this is indeed a norm on $\vemSpace \subset H^1_0(\Omega)$.

For the remainder of this paper, for ease of presentation we drop the dependence of $\mu$ on $\boldsymbol{x}$ and simply write $\mu(t)$ in place of $\mu(\boldsymbol{x},t)$.

\subsection{Approximation properties}\label{sec: approximation properties}
In order to carry out the a posteriori analysis, we require the following approximation result for the virtual element spaces defined in section~\ref{sec: virtual element method}.
The proof of the following result for the original VEM space of \cite{beirao_da_veiga_basic_2013} can be found in \cite{mora2015virtual} for the two dimensional case, and has been extended to three dimensions in \cite{cangiani2017posteriori} for the virtual element spaces of \cite{cangiani_conforming_2015} which are the same as those presented in this paper.
\begin{theorem}[Approximation using VEM functions]\label{thm: vem approximation}
    Under Assumption~\ref{ass: mesh regularity}, for any $w \in H^1(\Omega)$, there exists $w_I \in \vemSpace$ such that for all $\element \in \mesh$,
    \begin{align}\label{eqn: vem approximation}
        \| w - w_I \|_{0,\element} + h_{\element} | w- w_I |_{1,\element} \leq C_5 h_{\element} | w|_{1,\element},
    \end{align}
    where the constant $C_5$ depends only on $\polOrder$ and the mesh regularity.
\end{theorem}

\subsection{The residual equation}

We begin this section by deriving our \emph{residual equation}.
In order to do so, we define the \emph{error} $\error:= u-u_h \in H^1_0(\Omega)$, where $u$ is the weak solution to \eqref{eqn: var form pde} and $u_h$ is the virtual element solution to \eqref{eqn: vem problem}.
For any $\element \in \mesh$, using property assumption \eqref{eqn: mu prop 2}, we see that
\begin{align}\label{eqn: first estimate}
    C_2 |u - u_h|_{1}^2 &= C_2 \sum_{\element \in \mesh} \int_{\element} | \nabla (u-u_h) |^2 \, \dx
    \nonumber
    \\
    &\leq
    \sum_{\element \in \mesh} \int_{\element} \mu(|\nabla u|) \nabla u \cdot \nabla (u-u_h) - \mu (|\nabla u_h |) \nabla u_h  \cdot \nabla (u-u_h) \, \dx
    \nonumber
    \\
    &=
    \sum_{\element \in \mesh} \left[ a^{\element} (u;u,\error) - a^{\element}(u_h;u_h,\error) \right].
\end{align}
Therefore, since $u$ and $u_h$ are the solutions to \eqref{eqn: var form pde} and \eqref{eqn: vem problem}, respectively, for any $\chi \in \vemSpace$ we have
\begin{align}
    \sum_{\element \in \mesh} \big[ a^{\element} (u&;u,\error) - a^{\element}(u_h;u_h,\error)\big] \nonumber
    =
    \sum_{\element \in \mesh}\left[
    (f,\error)_{\element} - a^{\element}(u_h;u_h,\chi) - a^{\element}(u_h;u_h,\error-\chi)\right]
    \nonumber
    \\
    =& \
    \sum_{\element \in \mesh}
    \left[ (f,\error)_{\element} -(f_h,\chi)_{\element} +a_h^{\element}(u_h;u_h,\chi)- a^{\element}(u_h;u_h,\chi)
    - a^{\element}(u_h;u_h,\error-\chi) \right]
    \nonumber
    \\
    =& \
    \sum_{\element \in \mesh}
    \big[
    (f,\error- \chi)_{\element} + (f-f_h,\chi)_{\element}  +a_h^{\element}(u_h;u_h,\chi)
    \nonumber
    \\
    &\qquad\qquad-a^{\element}(u_h;u_h,\chi) - a^{\element}(u_h;u_h,\error-\chi)\big].
    \label{eqn: residual}
\end{align}

\subsection{Upper bounds (reliability)}
To derive computable (reliable) upper bounds on the error, we estimate each term in the residual equation \eqref{eqn: residual} in turn.
Now, for any $w \in H^1_0(\Omega)$, we introduce the gradient projection $\gradProj u_h$ of $u_h$ and recall that $\gradProj u_h = \cP^{\element}_{\polOrder-1} \nabla u_h$ (see \eqref{eqn: L2 property grad proj} in Lemma~\ref{lemma: l2 property}); hence,
\begin{align*}
    a^{\element}(u_h;u_h,w)
    =& \ \int_{\element} \mu(|\nabla u_h|) \nabla u_h \cdot \nabla w \, \dx
    \\
    =& \
    \int_{\element} \left( \mu(|\nabla u_h|) \nabla u_h \cdot \nabla w - \mu(| \gradProj u_h|) \gradProj u_h \cdot \nabla w \right) \, \dx
    \\
    &+ \int_{\element} \mu(| \gradProj u_h|) \gradProj u_h \cdot \nabla w \, \dx.
\intertext{Using integration by parts on the term on the second line, we see that}
    a^{\element}(u_h;u_h,w)
    =& \
    \int_{\element} ( \mu(|\nabla u_h|) \nabla u_h - \mu(|\cP^{\element}_{\polOrder-1} \nabla u_h |)  \cP^{\element}_{\polOrder-1} \nabla u_h )\cdot \nabla w \, \dx
    \\
    &- \int_{\element} \nabla \cdot (  \mu(| \cP^{\element}_{\polOrder-1} \nabla u_h |)  \cP^{\element}_{\polOrder-1} \nabla u_h ) w \, \dx
    \\
    &+ \sum_{e\subset \partial \element } \int_{e} \llbracket \mu(|\cP^{\element}_{\polOrder-1} \nabla u_h |) \cP^{\element}_{\polOrder-1} \nabla u_h \rrbracket w \, \ds.
\end{align*}
Furthermore, we introduce the polynomial approximation of order $\polOrder$, $\mu_h(t)=\cP^{\element}_{\polOrder} (\mu (t))$, of the coefficient $\mu$.
\begin{align*}
    a^{\element} (u_h;u_h&,w)
    =
    \int_{\element} ( \mu(|\nabla u_h|) \nabla u_h - \mu(|\cP^{\element}_{\polOrder-1} \nabla u_h|) \cP^{\element}_{\polOrder-1} \nabla u_h ) \cdot \nabla w \, \dx
    \\
    &+
    \int_{\element} \left( \nabla \cdot (\mu_h (|\cP^{\element}_{\polOrder-1}\nabla u_h | ) \cP^{\element}_{\polOrder-1} \nabla u_h )
    - \nabla \cdot ( \mu (| \cP^{\element}_{\polOrder-1} \nabla u_h | ) \cP^{\element}_{\polOrder-1} \nabla u_h  ) \right) w \, \dx
    \\
    &-
    \int_{\element} \nabla \cdot (\mu_h (|\cP^{\element}_{\polOrder-1} \nabla u_h | ) \cP^{\element}_{\polOrder-1} \nabla u_h ) w \, \dx
    \\
    &+
    \sum_{e \subset \partial \element } \Big( \int_e (\llbracket ( \mu(| \cP^{\element}_{\polOrder-1} \nabla u_h |) - \mu_h(|\cP^{\element}_{\polOrder-1}\nabla u_h |))\cP^{\element}_{\polOrder-1} \nabla u_h \rrbracket )  w \, \ds
    \\
    &+  \int_e (\llbracket \mu_h(|\cP^{\element}_{\polOrder-1}\nabla u_h |)\cP^{\element}_{\polOrder-1} \nabla u_h \rrbracket ) w \, \ds \Big).
\end{align*}

Substituting this into \eqref{eqn: residual} with $w=\error-\chi$, we see that
\begin{align}
    \sum_{\element \in \mesh}
    a^{\element} (u;u,\error) &- a^{\element}(u_h;u_h,\error)
    =
    \sum_{\element \in \mesh}
    \bigg(
    (f-f_h,\chi)_{\element} + a_h^{\element}(u_h;u_h,\chi) - a^{\element} (u_h;u_h,\chi)
    \nonumber
    \\
    &+
    (R^{\element},\error-\chi)_{\element} + (\theta^{\element},\error-\chi)_{\element} + (B^{\element}, \nabla (\error-\chi))_{\element}
    \nonumber
    \\
    &-
    \sum_{e \subset \partial \element} \left( (J^e,\error-\chi)_{0,e} + (\theta^e,\error-\chi)_{0,e} \right) \bigg)
    \label{eqn: final form residual}
\end{align}
where
\begin{align*}
    R^{\element} &:= (f_h + \nabla \cdot \mu_h(| \cP^{\element}_{\polOrder-1} \nabla u_h|) \cP^{\element}_{\polOrder-1} \nabla u_h )|_{\element},
    \\
    \theta^{\element} &:= ( f-f_h + \nabla \cdot (\mu(| \cP^{\element}_{\polOrder-1} \nabla u_h|)\cP^{\element}_{\polOrder-1} \nabla u_h - \mu_h (|\cP^{\element}_{\polOrder-1} \nabla u_h|) \cP^{\element}_{\polOrder-1} \nabla u_h))|_{\element},
    \\
    B^{\element} &:= (\mu(|\cP^{\element}_{\polOrder-1} \nabla u_h|) \cP^{\element}_{\polOrder-1} \nabla u_h - \mu(|\nabla u_h|) \nabla u_h)|_{\element},
    \\
    J^e &:= \llbracket \mu_h(| \cP^{\element}_{\polOrder-1} \nabla u_h |)\cP^{\element}_{\polOrder-1} \nabla u_h \rrbracket |_{e},
    \\
    \theta^e &:= \llbracket ( \mu(|\cP^{\element}_{\polOrder-1} \nabla u_h |) - \mu_h(|\cP^{\element}_{\polOrder-1} \nabla u_h |))\cP^{\element}_{\polOrder-1} \nabla u_h \rrbracket |_{e}.
\end{align*}

This leads us to the first crucial result of this subsection.
This proof follows the ideas in \cite{cangiani2017posteriori}.
\begin{theorem}[Upper bound]\label{thm: upper bound}
    Let $u\in H^1_0(\Omega)$ be the weak solution given by \eqref{eqn: var form pde} and $u_h \in \vemSpace$ be its virtual element approximation obtained from \eqref{eqn: vem problem}. Then, the following error bound holds:
    \begin{align*}
        | u - u_h|_{1}^2 \leq C \sum_{\element \in \mesh} (\eta_{\element}^2 + \Theta_{\element}^2 + \mathcal{S}_{\element}^2 + \Psi_{\element}^2)
    \end{align*}
    for some constant $C>0$ which is independent of $h$, $u$, and $u_h$, where
    \begin{align*}
        \eta_{\element}^2
        &:=
        h_{\element}^2 \| R^{\element} \|_{0,\element}^2 + \sum_{e \subset \partial \element} h_e \| J^e \|^2_{0,e},
        \\
        \Theta_{\element}^2
        &:=
        h_{\element}^2 \| \theta^{\element} \|_{0,\element}^2 + h_{\element}^2 \| f-f_h \|_{0,\element}^2 + \sum_{e \subset \partial \element} h_e \| \theta^e \|_{0,e}^2,
        \\
        \mathcal{S}_{\element}^2
        &:= S^{\element}(u_h; (I-\cP^{\element}_{\polOrder})u_h, (I-\cP^{\element}_{\polOrder})u_h),
        \\
        \Psi_{\element}^2
        &:=
        \| (\cP^{\element}_{\polOrder-1} -I) (\mu(|\cP^{\element}_{\polOrder-1} \nabla u_h|) \cP^{\element}_{\polOrder-1} \nabla u_h ) \|_{0,\element}^2 .
    \end{align*}
\end{theorem}

\begin{proof}
    Firstly, we take $\chi = \error_I \in \vemSpace$ in \eqref{eqn: final form residual} to be the interpolation of $\error=u-u_h$ into the VEM space $\vemSpace$; then, it follows from \eqref{eqn: first estimate} that
    \begin{align*}
        C_2 | \error |_{1}^2
        \leq \ &
        \sum_{\element \in \mesh } a^{\element}(u;u,\error) - a^{\element}(u_h;u_h,\error)
        \\
        = \ &
        \sum_{\element \in \mesh} \big[ (R^{\element},\error-\error_I)_{\element} + (\theta^{\element},\error-\error_I)_{\element} + (f-f_h,\error_I)_{\element} + (B^{\element},\nabla (\error-\error_I))_{\element}
        \\
        &+
        \left( a_h^{\element}(u_h;u_h,\error_I) - a^{\element}(u_h;u_h,\error_I) \right) \big]
        - \sum_{e \in \edgeSet} \left( (J^e ,\error-\error_I)_{0,e} + (\theta^e,\error-\error_I)_{0,e} \right)
        \intertext{Hence,}
        | \error |_{1}^2  \leq \ & C_2^{-1} \sum_{\element \in \mesh} (T_1^{\element} + T_2^{\element} + T_3^{\element} + T_4^{\element} + T_5^{\element} ) - C_2^{-1} \sum_{e \in \edgeSet} (T_6^e + T_7^e).
    \end{align*}
    We bound each term in turn, starting with $T_1^{\element}$ and $T_2^{\element}$.
    We use the interpolation approximation properties of VEM functions \eqref{eqn: vem approximation} detailed in Theorem~\ref{thm: vem approximation}, and Cauchy-Schwarz, which gives us
    \begin{align*}
        T_1^{\element} = (R^{\element},\error-\error_I)_{\element} \leq \| R^{\element} \|_{0,\element} \| \error-\error_I \|_{0,\element} &\leq C_5 h_{\element} \| R^{\element} \|_{0,\element} | \error |_{1,\element},
    \intertext{and}
        T_2^{\element} = (\theta^{\element},\error-\error_I)_{\element} \leq \| \theta^{\element} \|_{0,\element} \| \error- \error_I \|_{0,\element} &\leq C_5 h_{\element} \| \theta^{\element} \|_{\element} |\error|_{1,\element}.
    \end{align*}

    For $T_3^{\element}$, we use properties of the value projection and the $L^2$ projection from Theorem~\ref{thm: polynomial approximation}.
    Since $f_h := \cP^{\element}_{\polOrder} f$, we observe that $(f-f_h,\cP^{\element}_{\polOrder} \error_I) = 0$ and therefore
    \begin{align*}
        T_3^{\element} = (f-f_h,\error_I)_{\element}= (f-f_h,\error_I - \cP^{\element}_{\polOrder} \error_I) \leq \| f - f_h \|_{0,\element} C_3 h_{\element} | \error |_{1,\element}.
    \end{align*}

    Next, we look at bounding $T_4^{\element}$.
    For this term we use the property of $\mu$ detailed in \eqref{eqn: mu prop 1} as well as the interpolation results in \eqref{eqn: vem approximation}.
    Therefore
    \begin{align*}
        T_4^{\element} = (B^{\element} ,\nabla(\error-\error_I))_{\element}
        &= ( \mu(|\cP^{\element}_{\polOrder-1} \nabla u_h|) \cP^{\element}_{\polOrder-1} \nabla u_h - \mu(|\nabla u_h|) \nabla u_h, \nabla (\error - \error_I) )_{\element}
        \\
        &\leq \int_{\element} \big| \mu(|\cP^{\element}_{\polOrder-1} \nabla u_h|) \cP^{\element}_{\polOrder-1} \nabla u_h  - \mu(|\nabla u_h|) \nabla u_h \cdot \nabla (\error-\error_I) \big| \, \dx
        \\
        &\leq C_1 \| (\cP^{\element}_{\polOrder-1} - I) \nabla u_h \|_{0,\element} |\error-\error_I|_{1,\element}
        \\
        &\leq
        C_1 C_5   \| (\cP^{\element}_{\polOrder-1} - I) \nabla u_h \|_{0,\element} |\error|_{1,\element}.
    \end{align*}
    Since $\nabla \cP^{\element}_{\polOrder} u_h = \cP^{\element}_{\polOrder-1} \nabla \cP^{\element}_{\polOrder} u_h$, we use stability properties of the $L^2$ projection to see that
    \begin{align}
        \| (\cP^{\element}_{\polOrder-1}- I) \nabla u_h \|_{0,\element}
        &= \| (\cP^{\element}_{\polOrder-1} - I) \nabla u_h - (\cP^{\element}_{\polOrder-1} \nabla \cP^{\element}_{\polOrder} u_h - \nabla \cP^{\element}_{\polOrder} u_h ) \|_{0,\element}
        \nonumber
        \\
        &=
        \| (\cP^{\element}_{\polOrder-1}-I) \nabla ( I - \cP^{\element}_{\polOrder} ) u_h \|_{0,\element}
        \leq C
        \| \nabla ( I - \cP^{\element}_{\polOrder} ) u_h  \|_{0,\element}.
        \label{eqn: bound before stab}
    \end{align}
    However, since $\nabla u_h$ is not a computable quantity, we need to bound this in terms of the stabilisation $S^{\element}$.
    We apply \eqref{eqn: stab bound} with $z_h = u_h$ and $v_h = (I-\cP^{\element}_{\polOrder})u_h$ to see that
    \begin{align}\label{eqn: bound by stabilisation}
        \| \nabla (I-\cP^{\element}_{\polOrder})u_h \|_{0,\element} \leq C_4 (S^{\element}(u_h; (I-\cP^{\element}_{\polOrder})u_h, (I-\cP^{\element}_{\polOrder})u_h))^{\half}.
    \end{align}
    Therefore,
    \begin{align*}
        T_4^{\element} \leq C_1 C_4 C_5 |\error|_{1,\element} (S^{\element}(u_h; (I-\cP^{\element}_{\polOrder})u_h, (I-\cP^{\element}_{\polOrder})u_h))^{\half}.
    \end{align*}
    Before we can estimate the edge terms, we turn our attention to $T_5^{\element}$.
    We use property \eqref{eqn: mu prop 1}, together with properties of the $L^2$ projection to see that
    \begin{align*}
        T_5^{\element} :=& \ a_h^{\element}(u_h;u_h,\error_I) - a^{\element}(u_h;u_h,\error_I)
        \\
        =& \ \int_{\element}  \mu(|\cP^{\element}_{\polOrder-1} \nabla u_h|) \cP^{\element}_{\polOrder-1} \nabla u_h \cdot \cP^{\element}_{\polOrder-1} \nabla \error_I - \mu(|\nabla u_h |) \nabla u_h \cdot \nabla \error_I  \, \dx
        \\
        &+ S^{\element} ( u_h; (I-\valueProj) u_h, (I-\valueProj)\error_I)
        \\
        =& \
        \int_{\element} (\mu(|\cP^{\element}_{\polOrder-1} \nabla u_h|) \cP^{\element}_{\polOrder-1} \nabla u_h -  \mu(|\nabla u_h |) \nabla u_h ) \cdot  \nabla \error_I \, \dx
        \\
        &+
        \int_{\element} \mu(|\cP^{\element}_{\polOrder-1} \nabla u_h|) \cP^{\element}_{\polOrder-1} \nabla u_h \cdot (\cP^{\element}_{\polOrder-1} - I) \nabla \error_I  \, \dx
        + S^{\element} ( u_h; (I-\valueProj) u_h, (I-\valueProj)\error_I)
        \\
        \leq& \
        \left( C_1 \| \cP^{\element}_{\polOrder-1} \nabla u_h - \nabla u_h\|_{0,\element} + \| (\cP^{\element}_{\polOrder-1} - I)\mu(|\cP^{\element}_{\polOrder-1} \nabla u_h|) \cP^{\element}_{\polOrder-1} \nabla u_h \|_{0,\element} \right) |\error_I |_{1,\element}
        \\
        &+ S^{\element} ( u_h; (I-\valueProj) u_h, (I-\valueProj)\error_I).
    \end{align*}
    Using Cauchy-Schwarz, we see that
    \begin{align*}
        S^{\element} ( u_h; (I-\valueProj) u_h, (I-\valueProj)\error_I)
        \leq& \
        (S^{\element} ( u_h; (I-\valueProj) u_h, (I-\valueProj) u_h))^{\half}
        \\
        &\times (S^{\element} ( u_h; (I-\valueProj) \error_I, (I-\valueProj) \error_I ))^{\half}.
    \end{align*}
    In order to bound the second part of this term, we use \eqref{eqn: mu bound}. Then,
    \begin{align*}
        S^{\element} ( u_h; (I-\valueProj) \error_I, (I-\valueProj) \error_I )
        &\leq
        C^* \int_{\element} \mu(|\nabla u_h|) \nabla (I-\valueProj) \error_I \cdot \nabla (I-\valueProj) \error_I \, \dx
        \\
        &\leq
        C^* M_{\mu} \| \nabla (I-\valueProj) \error_I \|_{0,\element}^2
        \\
        &\leq C^* M_{\mu} ( C_3 |\error_I|_{1,\element} )^2
    \end{align*}
    where we have used Theorem~\ref{thm: polynomial approximation} as well as the stability property detailed in \eqref{eqn: admissible stabilisation} in the first line.
    Therefore,
    \begin{align*}
        S^{\element} ( u_h; (I-\valueProj) u_h, (I-\valueProj)\error_I)
        \leq
        (C^* M_{\mu})^{\half} C_3 |\error_I|_{1,\element}
        (S^{\element} ( u_h; (I-\valueProj) u_h, (I-\valueProj) u_h))^{\half}.
    \end{align*}
    Recalling \eqref{eqn: bound by stabilisation}, we have
    \begin{align*}
        T_5^{\element}
        \leq& \
        |\error_I |_{1,\element}
        \left( C_1 \| \cP^{\element}_{\polOrder-1} \nabla u_h - \nabla u_h\|_{0,\element} + \| (\cP^{\element}_{\polOrder-1} - I)\mu(|\cP^{\element}_{\polOrder-1} \nabla u_h|) \cP^{\element}_{\polOrder-1} \nabla u_h \|_{0,\element} \right.
        \\
        &+
        \left.
        (C^* M_{\mu})^{\half} C_3 (S^{\element} ( u_h; (I-\valueProj) u_h, (I-\valueProj) u_h))^{\half} \right)
        \\
        \leq& \
        |\error_I|_{1,\element} (
            C_1 C_4 (S^{\element}(u_h; (I-\cP^{\element}_{\polOrder})u_h, (I-\cP^{\element}_{\polOrder})u_h))^{\half}
        \\
        &+
        \| (\cP^{\element}_{\polOrder-1}-I )\mu(|\cP^{\element}_{\polOrder-1} \nabla u_h|) \cP^{\element}_{\polOrder-1} \nabla u_h \|_{0,\element}
        \\
        &+
        (C^* M_{\mu})^{\half} C_3  (S^{\element} ( u_h; (I-\cP^{\element}_{\polOrder}) u_h, (I-\cP^{\element}_{\polOrder}) u_h))^{\half}
        )
    \end{align*}
    where we have used the $L^2$ property of the value projection \eqref{eqn: L2 property value proj} in the last step.
    Hence,
    \begin{align*}
        T_5^{\element} \leq C |\error_I|_{1,\element} \left( ( S^{\element} ( u_h; (I-\cP^{\element}_{\polOrder}) u_h, (I-\cP^{\element}_{\polOrder}) u_h))^{\half} \right.
        \\
        \left.
        + \| (\cP^{\element}_{\polOrder-1}-I )\mu(|\cP^{\element}_{\polOrder-1} \nabla u_h|) \cP^{\element}_{\polOrder-1} \nabla u_h \|_{0,\element} \right).
    \end{align*}

    Lastly, we bound $T_6^e$ and $T_7^e$.
    This requires the use of the standard scaled trace inequality which states that: for $v \in H^1(\element)$
    \begin{align*}
        \| v \|^2_{0,e} \leq C_6 ( h_{\element}^{-1} \| v\|^2_{0,\element} + h_{\element} \| \nabla v \|^2_{0,\element}).
    \end{align*}
    This, together with \eqref{eqn: vem approximation}, gives us
    \begin{align*}
        T_6^e = \int_{e} J^e (\error-\error_I) \, \ds \leq \| J^e \|_{0,e} \| \error- \error_I\|_{0,e}
        \leq
        C_6^{\half} C_5 h_e^{\half} \rho^{-\half} |\error|_{1,E^{+} \cup E^{-}} \| J^e\|_{0,e}
        \intertext{and}
        T_7^e := \int_e \theta^e (\error-\error_I) \, \ds \leq \| \theta^{e} \|_{0,e} \| \error-\error_I\|_{0,e}
        \leq
        C_6^{\half} C_5 h_e^{\half} \rho^{-\half} |\error|_{1,E^{+} \cup E^{-}} \| \theta^e\|_{0,e}
    \end{align*}
    where we have also applied part \ref{ass: mesh reg part b} of the mesh regularity Assumption~\ref{ass: mesh regularity}.

    The result now follows from Young's inequality since we have that
    \begin{align*}
        |u-u_h|_1 \leq& \
        C \Bigg\{ \left( \sum_{\element \in \mesh} h_{\element}^2 \| R^{\element}\|_{0,\element}^2 \right)^{\half}
        +
        \left( \sum_{e \in \edgeSet} h_{e} \| J^{e}\|_{0,e}^2 \right)^{\half}
        \\
        &+
        \left( \sum_{\element \in \mesh} h_{\element}^2 \| \theta^{\element}\|_{0,\element}^2 \right)^{\half}
        +
        \left( \sum_{\element \in \mesh} h_{\element}^2 \| f-f_h\|_{0,\element}^2 \right)^{\half}
        +
        \left( \sum_{e \in \edgeSet} h_{e} \| \theta^{e}\|_{0,e}^2 \right)^{\half}
        \\
        &+
        \left( \sum_{\element \in \mesh} S^{\element}(u_h; (I-\cP^{\element}_{\polOrder})u_h, (I-\cP^{\element}_{\polOrder})u_h)\right)^{\half}
        \\
        &+ \left( \sum_{\element \in \mesh} \| (\cP^{\element}_{\polOrder-1} -I)\mu(|\cP^{\element}_{\polOrder-1} \nabla u_h|) \cP^{\element}_{\polOrder-1} \nabla u_h \|_{0,\element}^2  \right)^{\half} \Bigg\}.
    \end{align*}
\end{proof}

\begin{remark}
    We note that the estimator in Theorem~\ref{thm: upper bound} is also an estimator for the projected solution, as detailed in the next corollary.
\end{remark}

\begin{corollary}\label{cor: upper bound proj}
    Let $u\in H^1_0(\Omega)$ be the weak solution given by \eqref{eqn: var form pde} and $u_h \in \vemSpace$ be its virtual element approximation obtained from \eqref{eqn: vem problem}. Then, the following error bounds hold:
    \begin{align}
        | u - \valueGlobal u_h |_{1}^2 &\leq \overline{C} \sum_{\element \in \mesh} (\eta_{\element}^2 + \Theta_{\element}^2 + \mathcal{S}_{\element}^2 + \Psi_{\element}^2)
        \label{eqn: corollary for value projection}
        \\
        \| \nabla u - \Pi^h_1 u_h \|_{0}^2 &\leq \widehat{C} \sum_{\element \in \mesh} (\eta_{\element}^2 + \Theta_{\element}^2 + \mathcal{S}_{\element}^2 + \Psi_{\element}^2)
        \label{eqn: corollary for grad projection}
    \end{align}
    for constants $\overline{C},\widehat{C}>0$ which are independent of $h$, $u$, and $u_h$, where $\eta_{\element}, \Theta_{\element}, \mathcal{S}_{\element}$, and $\Psi_{\element}$ are defined in Theorem~\ref{thm: upper bound}.
\end{corollary}

\begin{proof}
    First, we show \eqref{eqn: corollary for value projection} using \eqref{eqn: stab bound} with $z_h = u_h$, and $v_h = u_h - \cP^{\element}_{\polOrder} u_h$.
    Together with $L^2$ projection properties of the value projection, we notice that
    \begin{align*}
        | u_h - \valueGlobal u_h |^2_{1}
        &=
        \sum_{\element \in \mesh} \| \nabla (u_h - \cP^{\element}_{\polOrder} u_h ) \|^2_{0,\element}
        \\
        &\leq
        \sum_{\element \in \mesh}
        C_4 S^{\element}(u_h;u_h - \cP^{\element}_{\polOrder} u_h ,u_h - \cP^{\element}_{\polOrder} u_h).
    \end{align*}
    Therefore, applying the triangle inequality, we see that
    \begin{align*}
        | u - \valueGlobal u_h |^2_{1}
        &\leq
        2 | u - u_h |^2_{1} + 2| u_h - \valueGlobal u_h |^2_{1}
        \\
        &\leq
        2|u-u_h|^2_{1} + 2 \sum_{\element \in \mesh}
        C_4 S^{\element}(u_h;u_h -\cP^{\element}_{\polOrder}  u_h,u_h - \cP^{\element}_{\polOrder}  u_h).
    \end{align*}
    The result now follows from Theorem~\ref{thm: upper bound}.

    To show \eqref{eqn: corollary for grad projection}, we note that
    \begin{align*}
        \| \nabla u_h - \Pi^h_1 u_h \|_0^2 = \sum_{\element \in \mesh} \| \nabla u_h - \gradProj u_h \|_{0,\element}^2
        = \sum_{\element \in \mesh} \| (I  - \cP_{\polOrder-1}^{\element} ) \nabla u_h \|_{0,\element}^2
    \end{align*}
    where we have applied Lemma~\ref{lemma: l2 property}.
    Therefore, applying \eqref{eqn: bound before stab} and \eqref{eqn: bound by stabilisation} gives us the desired result as before.
\end{proof}

\subsection{Lower bounds (efficiency)}
In this subsection we prove local lower bounds of the error in the $H^1$-seminorm which demonstrates the efficiency of the error bound.
In order to prove these lower bounds, we use standard properties of bubble functions introduced below.

A bubble function $\psi^{\element} \in H^1_0(\element)$ is constructed piecewise as the sum of the bubble functions \cite{ainsworth1997posteriori} on each triangle of the mesh sub-triangulation \cite{cangiani2017posteriori}; cf. Remark~\ref{rmk: subtriangulation}.

\begin{lemma}[Bubble functions]\label{lemma: interior bubbles}
    For $\element \in \mesh$, let $\psi^{\element}$ be the corresponding bubble function.
    Then, there exists a constant $C_7$, independent of $h_{\element}$ such that for all $p \in \prob_{\polOrder}(\element)$
    \begin{align*}
        C_7^{-1} \| p \|^2_{0,\element} \leq \int_{\element} \psi^{\element} p^2 \, \dx \leq C_7 \| p \|^2_{0,\element},
    \end{align*}
    and
    \begin{align*}
        C_7^{-1} \| p \|_{0,\element} \leq \| \psi^{\element} p \|_{0,\element} + h_{\element} | \psi^{\element} p |_{1,\element} \leq C_7 \| p \|_{0,\element}.
    \end{align*}

    Furthermore, for each $e \subset \partial \element$, let $\psi^e$ denote the corresponding edge bubble function.
    Then, for all $p \in \prob_{\polOrder}(e)$,
    \begin{align*}
        C_7^{-1} \| p \|^2_{0,e} \leq \int_{e} \psi^{e} p^2 \, \ds \leq C_7 \| p \|^2_{0,e},
    \end{align*}
    and
    \begin{align*}
        h_{\element}^{-\half} \| \psi^{e} p \|_{0,\element} + h_{\element}^{\half} | \psi^{e} p |_{1,\element} \leq C_7 \| p \|_{0,e}
    \end{align*}
    for a constant $C_7$ independent of $h_{\element}$.
\end{lemma}

\begin{theorem}[Local lower bound]\label{thm: lower bound}
    Let $\eta_{\element}, \mathcal{S}_{\element},$ and $\Theta_{\element}$ be defined as in Theorem~\ref{thm: upper bound}.
    Then, there exists a constant $C>0$, independent of $h,u,$ and $u_h$, such that for $E\in\mesh$
    \begin{align*}
        \eta_{\element}^2 \leq C \sum_{E^{\prime} \in \omega_{\element}} \left( \| \nabla (u-u_h) \|^2_{0,\element^{\prime}} + \mathcal{S}_{\element^{\prime}}^2 + \Theta_{\element^{\prime}}^2 \right)
    \end{align*}
    where $\omega_{\element}$ denotes the patch of elements containing $\element$ and its neighbouring elements; i.e.,
    $$\omega_{\element} := \{ \element^{\prime} \in \mesh : \partial \element^{\prime} \cap \partial \element \neq \emptyset \}.$$
\end{theorem}
\begin{proof}
    Since $u \in H^1_0(\Omega)$ solves \eqref{eqn: var form pde}, it follows that
    \begin{align*}
        \sum_{\element \in \mesh} a^{\element} (u;u,\psi^{\element} R^{\element}) &- a^{\element}(u_h;u_h,\psi^{\element} R^{\element}) \nonumber
        =
        \sum_{\element \in \mesh}
        (f,\psi^{\element} R^{\element})_{\element} - a^{\element}(u_h;u_h,\psi^{\element} R^{\element}).
    \end{align*}
    We can now use the same approach that was used to derive \eqref{eqn: final form residual} but with $w=\psi^{\element} R^{\element}$ to see that
    \begin{align*}
        a^{\element} (u_h;u_h&,\psi^{\element} R^{\element})
        =
        \int_{\element} ( \mu(|\nabla u_h|) \nabla u_h - \mu(|\cP^{\element}_{\polOrder-1} \nabla u_h|) \cP^{\element}_{\polOrder-1} \nabla u_h ) \cdot \nabla (\psi^{\element} R^{\element}) \, \dx
        \\
        &+
        \int_{\element} \left( \nabla \cdot (\mu_h (|\cP^{\element}_{\polOrder-1}\nabla u_h | ) \cP^{\element}_{\polOrder-1} \nabla u_h )
        - \nabla \cdot ( \mu (| \cP^{\element}_{\polOrder-1} \nabla u_h | ) \cP^{\element}_{\polOrder-1} \nabla u_h  ) \right) \psi^{\element} R^{\element} \, \dx
        \\
        &-
        \int_{\element} \nabla \cdot (\mu_h (|\cP^{\element}_{\polOrder-1} \nabla u_h | ) \cP^{\element}_{\polOrder-1} \nabla u_h ) \psi^{\element} R^{\element} \, \dx
    \end{align*}
    where the boundary terms that appear in \eqref{eqn: final form residual} are zero since $\psi^{\element} |_{\partial \element} = 0$. Therefore,
    \begin{align*}
        a^{\element} (u;u,\psi^{\element} R^{\element}) - a^{\element}(u_h;u_h,\psi^{\element} R^{\element})
        =& \
        (R^{\element}, \psi^{\element} R^{\element})_{\element} + (\theta^{\element},\psi^{\element} R^{\element})_{\element} + (B^{\element}, \nabla (\psi^{\element} R^{\element}))_{\element}.
    \end{align*}
    We first observe using Assumption~\ref{ass: mu} that it also holds that
    \begin{align*}
        a^{\element} (u;u,\psi^{\element} R^{\element}) - a^{\element}(u_h;u_h,\psi^{\element} R^{\element})
        &\leq
        \int_{\element} \left| \left( \mu(|\nabla u|) \nabla u - \mu(|\nabla u_h|)\nabla u_h \right) \cdot \nabla (\psi^{\element} R^{\element}) \right| \, \dx
        \\
        &\leq
        C_1 \| \nabla (u - u_h) \|_{0,\element} \| \nabla (\psi^{\element} R^{\element}) \|_{0,\element}.
    \end{align*}

    Using the same argument as $T_4^{\element}$ in Theorem~\ref{thm: upper bound}, we can show that
    \begin{align*}
        (B^{\element},\nabla(\psi^{\element}R^{\element}))_{\element}
        &\leq
        C_1 \| (I-\cP^{\element}_{\polOrder-1})\nabla u_h \|_{0,\element} | \psi^{\element}R^{\element} |_{1,\element}
        \\
        &\leq
        C_1 C_4 (S^{\element}(u_h; (I-\cP^{\element}_{\polOrder})u_h, (I-\cP^{\element}_{\polOrder})u_h))^{\half} \| \nabla (\psi^{\element}R^{\element}) \|_{0,\element}.
    \end{align*}
    Next, using Lemma~\ref{lemma: interior bubbles} with $p=R^{\element} \in \prob_{k}(\element)$ for some $k \in \N$, we have that
    \begin{align*}
        C_7^{-1} \| R^{\element} \|_{0,\element}^2
        \leq& \
        (R^{\element},\psi^{\element}R^{\element})_{\element}
        \\
        =& \
        a^{\element} (u;u,\psi^{\element} R^{\element}) - a^{\element}(u_h;u_h,\psi^{\element} R^{\element})
        - (\theta^{\element},\psi^{\element} R^{\element})_{\element}
        - (B^{\element}, \nabla (\psi^{\element} R^{\element}))_{\element}
        \\
        \leq& \
        C_1 \| \nabla (u-u_h) \|_{0,\element} \| \nabla (\psi^{\element}R^{\element}) \|_{0,\element}
        +
        \| \theta^{\element} \|_{0,\element} \| \psi^{\element} R^{\element} \|_{0,\element}
        \\
        &+
        C_1 C_4 (S^{\element}(u_h; (I-\cP^{\element}_{\polOrder})u_h, (I-\cP^{\element}_{\polOrder})u_h))^{\half} \| \nabla (\psi^{\element}R^{\element}) \|_{0,\element}
        \\
        \leq& \
        \bigg( C_1 C_7 h_{\element}^{-1} \| \nabla \error \|_{0,\element} + C_7 \| \theta^{\element} \|_{0,\element}
        \\
        &+
        C_1 C_4 C_7 h_{\element}^{-1}  (S^{\element}(u_h; (I-\cP^{\element}_{\polOrder})u_h, (I-\cP^{\element}_{\polOrder})u_h))^{\half} \bigg) \| R^{\element} \|_{0,\element}
    \end{align*}
    where we have used Lemma~\ref{lemma: interior bubbles} again in the last line.
    Therefore
    \begin{align}\label{eqn: element residual}
        \tilde C h_{\element}^2 \| R^{\element} \|_{0,\element}^2 \leq  \| \nabla \error \|_{0,\element}^2 + h_{\element}^2 \| \theta^{\element} \|_{0,\element}^2 + \mathcal{S}_{\element}^2
    \end{align}
    for some constant $\tilde{C}$.

    For the edge residual part of $\eta^2_{\element}$, we extend $J^e$ into $\omega_e$ through a constant prolongation in the normal direction to the edge $e$; see e.g. \cite{cangiani2017posteriori}.
    This gives us $J^e \in \prob_{\polOrder}(\omega_e) \subset \vemSpace^{\omega_e} := \vemSpace^{\element^{+}} \cap \vemSpace^{\element^{-}}$ where $E^{+} \cap E^{-} = e$.
    Then, using the derivation of \eqref{eqn: final form residual} again with $w=\psi^e J^e$ and noting that the edge terms are no longer zero,
    \begin{align*}
        \sum_{\element \in \mesh} a^{\element}(u;u,\psi^e J^e) -& \ a^{\element}(u_h;u_h,\psi^e J^e )
        \\
        =& \
        \sum_{\element^{\prime} \in \omega_e} (R^{\element^{\prime}},\psi^e J^e)_{\element^{\prime}} + (\theta^{\element^{\prime}},\psi^e J^e)_{\element^{\prime}} + (B^{\element^{\prime}},\nabla(\psi^e J^e))_{\element^{\prime}}
        \\
        &-
        (J^e,\psi^e J^e)_{0,e} - (\theta^e, \psi^e J^e)_{0,e}.
    \end{align*}
    We use the same argument as before for the inner residual terms, as well as Lemma~\ref{lemma: interior bubbles}, to see that
    \begin{align*}
        C_7^{-1} \| J^e \|^2_{0,e} \leq& \ (J^e,\psi^e J^e)_{0,e}
        \\
        =& \
        \sum_{\element^{\prime} \in \omega_e} (R^{\element^{\prime}},\psi^e J^e)_{\element^{\prime}} + (\theta^{\element^{\prime}},\psi^e J^e)_{\element^{\prime}} + (B^{\element^{\prime}},\nabla(\psi^e J^e))_{\element^{\prime}}
        \\
        &-
        (\theta^e,\psi^e J^e)_{0,e} +  a(u_h;u_h,\psi^e J^e ) - a(u;u,\psi^e J^e)
        \\
        \leq& \
        \sum_{\element^{\prime} \in \omega_e} \bigg( ( \| R^{\element^{\prime}} \|_{0,\element^{\prime}} + \| \theta^{\element^{\prime}} \|_{0,\element^{\prime}} ) \| \psi^e J^e \|_{0,\element^{\prime}}
        \\
        &+ \left( C_1 C_4 \mathcal{S}_{\element^{\prime}} + C_1 \| \nabla \error \|_{0,\element^{\prime}} \right) \| \nabla (\psi^e J^e) \|_{0,\element^{\prime}} \bigg)
        + \| \theta^e \|_{0,e} \| \psi^e J^e \|_{0,e}
        \\
        \leq& \
        \sum_{\element^{\prime} \in \omega_e} \bigg(
        h_{\element^{\prime}}^{\half} ( \| R^{\element^{\prime}} \|_{0,\element^{\prime}} + \| \theta^{\element^{\prime}} \|_{0,\element^{\prime}} ) \| J^e \|_{0,e}
        \\
        &+
        h_{\element^{\prime}}^{-\half}
        \left( C_1 C_4 \mathcal{S}_{\element^{\prime}} + C_1 \| \nabla \error \|_{0,\element^{\prime}} \right) \|  J^e \|_{0,e} \bigg)
        + \| \theta^e \|_{0,e} \| J^e \|_{0,e}.
    \end{align*}
    Firstly, we divide through by $\| J^e \|_{0,e}$ to see that
    \begin{align*}
        C_7^{-1} \| J^e \|_{0,e}
        \leq
        \| \theta^e \|_{0,e}  +
        \sum_{\element^{\prime} \in \omega_e}
        h_{\element^{\prime}}^{\half} ( \| R^{\element^{\prime}} \|_{0,\element^{\prime}} + \| \theta^{\element^{\prime}} \|_{0,\element^{\prime}} )
        +
        h_{\element^{\prime}}^{-\half}
        \left( C_1 C_4 \mathcal{S}_{\element^{\prime}} + C_1 \| \nabla \error \|_{0,\element^{\prime}} \right).
    \end{align*}
    Therefore, applying \eqref{eqn: element residual}, we get the following
    \begin{align*}
        \tilde{C} \|  J^e \|^2_{0,e}
        \leq& \
        \| \theta^e \|_{0,e}^2  + \sum_{\element^{\prime} \in \omega_e}
        h_{\element^{\prime}} ( \| R^{\element^{\prime}} \|^2_{0,\element^{\prime}} + \| \theta^{\element^{\prime}} \|^2_{0,\element^{\prime}} )
        +
        h_{\element^{\prime}}^{-1}
        \left( \mathcal{S}_{\element^{\prime}}^2 + \| \nabla \error \|^2_{0,\element^{\prime}} \right)
        \\
        \leq& \  \| \theta^e \|_{0,e}^2  + \sum_{\element^{\prime} \in \omega_e}
        h_{\element^{\prime}}  \| \theta^{\element^{\prime}} \|^2_{0,\element^{\prime}}
        +
        2 h_{\element^{\prime}}^{-1}
        \left( \mathcal{S}_{\element^{\prime}}^2 +  \| \nabla \error \|^2_{0,\element^{\prime}} \right)
    \end{align*}
    for some constant $\tilde{C}$.
    Hence, multiplying by $h_e$, and observing that $h_e \leq h_{\element^{\prime}}$, we see that
    \begin{align*}
        \tilde{C} h_e \|  J^e \|^2_{0,e}
        \leq& \
        h_e \| \theta^e \|_{0,e}^2  + \sum_{\element^{\prime} \in \omega_e}
        h_e h_{\element^{\prime}}  \| \theta^{\element^{\prime}} \|^2_{0,\element^{\prime}}
        +
        2 h_e h_{\element^{\prime}}^{-1}
        \left( \mathcal{S}_{\element^{\prime}}^2 + \| \nabla \error \|^2_{0,\element^{\prime}} \right)
        \\
        =& \
        h_e \| \theta^e \|_{0,e}^2  + \sum_{\element^{\prime} \in \omega_e}
        h_{\element^{\prime}}^2 \| \theta^{\element^{\prime}} \|^2_{0,\element^{\prime}}
        +
        2 ( \mathcal{S}_{\element^{\prime}}^2 +  \| \nabla \error \|^2_{0,\element^{\prime}} )
    \end{align*}
    Recalling the definition of $\eta_{\element}^2$ and combining the estimates gives us the desired result.
\end{proof}

\begin{remark}
    We note that a direct consequence of Theorem~\ref{thm: lower bound} is a lower bound on the error between the solution $u$ and the projected solution. That is,
    \begin{align}\label{eqn: corollary lower bound value proj}
        \eta_{\element}^2 \leq C \sum_{E^{\prime} \in \omega_{\element}} \left( | u - \Pi_0^{\element^{\prime}}u_h |^2_{1,\element^{\prime}} + \mathcal{S}_{\element^{\prime}}^2 + \Theta_{\element^{\prime}}^2 \right)
    \intertext{and}
        \eta_{\element}^2 \leq C \sum_{E^{\prime} \in \omega_{\element}} \left( \| \nabla u - \Pi_1^{\element^{\prime}} u_h \|^2_{0,\element^{\prime}} + \mathcal{S}_{\element^{\prime}}^2 + \Theta_{\element^{\prime}}^2 \right)
        \label{eqn: corollary lower bound grad proj}
    \end{align}
    with $\omega_{\element}$ defined as in Theorem~\ref{thm: lower bound}.
    Note that \eqref{eqn: corollary lower bound value proj} and \eqref{eqn: corollary lower bound grad proj} follow directly from Theorem~\ref{thm: lower bound} combined with an application of the triangle inequality.
\end{remark}

In the last part of this subsection, we prove the following lower bound for the inconsistency term $\Psi_{\element}^2$, where $\Psi_{\element}^2$ is defined in Theorem~\ref{thm: upper bound}.

\begin{theorem}[Lower bound for inconsistency term]
    There exists a constant $C>0$ independent of $h,u,$ and $u_h$ such that
    \begin{align*}
        \Psi_{\element}^2 \leq C \left( \| \nabla( u - u_h ) |^2_{0,\element} + \mathcal{S}_{\element}^2 + \| ( \cP^{\element}_{\polOrder-1} - I) \mu(|\nabla u|) \nabla u \|^2_{0,\element} \right).
    \end{align*}
\end{theorem}
\begin{proof}
    Recall the definition of $\Psi_{\element}^2$ (Theorem~\ref{thm: upper bound}), then note that we have
    \begin{align*}
        \Psi_{\element}^2
        =& \
        \| (\cP^{\element}_{\polOrder-1} -I) (\mu(|\cP^{\element}_{\polOrder-1} \nabla u_h|) \cP^{\element}_{\polOrder-1} \nabla u_h ) \|_{0,\element}^2
        \\
        =& \
        \big( \| \cP^{\element}_{\polOrder-1} (\mu(|\cP^{\element}_{\polOrder-1} \nabla u_h|)  \cP^{\element}_{\polOrder-1} \nabla u_h ) -
        \cP^{\element}_{\polOrder-1} ( \mu(|\nabla u|) \nabla u ) \|_{0,\element}
        \\
        &+ \| \cP^{\element}_{\polOrder-1} ( \mu(|\nabla u|) \nabla u ) -  \mu(| \nabla u|)  \nabla u \|_{0,\element}
        \\
        &+
        \| \mu(|\nabla u |) \nabla u -  \mu(|\cP^{\element}_{\polOrder-1} \nabla u_h|) \cP^{\element}_{\polOrder-1}  \nabla u_h  \|_{0,\element} \big)^2
        \\
        \leq& \
        2  \| \cP^{\element}_{\polOrder-1} (\mu(|\cP^{\element}_{\polOrder-1} \nabla u_h|)  \cP^{\element}_{\polOrder-1} \nabla u_h ) -
        \cP^{\element}_{\polOrder-1} ( \mu(|\nabla u|) \nabla u ) \|_{0,\element}^2
        \\
        &+
        2\| \cP^{\element}_{\polOrder-1} ( \mu(|\nabla u|) \nabla u ) -  \mu(| \nabla u|)  \nabla u \|_{0,\element}^2
        \\
        &+
        2 \| \mu(|\nabla u |) \nabla u -  \mu(|\cP^{\element}_{\polOrder-1} \nabla u_h|) \cP^{\element}_{\polOrder-1}  \nabla u_h  \|_{0,\element}^2
        \\
        \leq& \
        2 \| ( \cP^{\element}_{\polOrder-1} - I) \mu(|\nabla u|) \nabla u \|^2_{0,\element}
        +
        4 \| \mu(|\nabla u |) \nabla u -  \mu(|\cP^{\element}_{\polOrder-1} \nabla u_h|) \cP^{\element}_{\polOrder-1}  \nabla u_h  \|_{0,\element}^2
    \end{align*}
    where we have used stability properties of the $L^2$ projection in the last step.
    We bound the last term by first applying \eqref{eqn: mu prop 1} to see that
    \begin{align*}
        \| \mu(|\nabla u |) \nabla u &-  \mu(|\cP^{\element}_{\polOrder-1} \nabla u_h|) \cP^{\element}_{\polOrder-1}  \nabla u_h  \|_{0,\element}^2
        \leq C_1 \|  \nabla u - \cP^{\element}_{\polOrder-1}  \nabla u_h \|^2_{0,\element}
        \\
        &\leq 2 C_1 ( \| \nabla u - \nabla u_h \|_{0,\element}^2 + \| \nabla u_h - \cP^{\element}_{\polOrder-1}  \nabla u_h \|_{0,\element}^2).
    \end{align*}
    Using the same argument for term $T_4^{\element}$ in Theorem~\ref{thm: upper bound}, see e.g. \eqref{eqn: bound by stabilisation}, we can show that
    \begin{align*}
        \| (I-\cP^{\element}_{\polOrder-1} ) \nabla u_h \|_{0,\element} \leq C_4 (S^{\element}(u_h; (I-\cP^{\element}_{\polOrder})u_h, (I-\cP^{\element}_{\polOrder})u_h))^{\half}.
    \end{align*}
    Therefore, the result follows.
\end{proof}

\section{Numerical results}\label{sec: numerics}

In this section we present a collection of numerical results aimed at investigating the behaviour of the a posteriori error estimator derived in Theorem~\ref{thm: upper bound} and Corollary~\ref{cor: upper bound proj}.
Furthermore, we present an estimator driven adaptive algorithm and apply it to a selection of test problems.

The code used to carry out the simulations is based on the Distributed and Unified Numerics Environment (DUNE) software framework \cite{dunegridpaperII}.
The virtual element method has been implemented within the DUNE-FEM module \cite{dedner2010generic} and further implementation details can be found in \cite{dedner2022framework}.
DUNE is open source software implemented in C++; however, a user can readily carry out numerical experiments by describing mathematical models using the domain specific form language UFL \cite{alnaes_unified_2012} within the Python frontend \cite{dedner2020python,dedner_dune_2018}.

The aim of these experiments is to demonstrate that the second a posteriori error indicator \eqref{eqn: corollary for grad projection} in Corollary~\ref{cor: upper bound proj} converges at the same rate as $\| \nabla u - \gradProj u_h \|_{0,\Omega}$ on a sequence of adaptively refined meshes.
We note that we consider this bound rather than Theorem~\ref{thm: upper bound} as $\nabla u_h$ is \emph{not} a computable quantity. As is standard in residual-based a posteriori error estimation we set the constant $\widehat{C}$ in Corollary~\ref{cor: upper bound proj} to $1$ for simplicity; in general, this constant should be determined numerically, cf. \cite{eriksson1995}. We are then able to check whether the effectivity index,
\begin{align}
    \text{effectivity} := \frac{\left( \sum_{\element \in \mesh} \eta^2_{\element} + \Theta^2_{\element} + \mathcal{S}^2_{\element} + \Psi^2_{\element} \right)^{\half}}{ \| \nabla u - \gradProj u_h \|_{0,\Omega} }
    \label{eqn: effectivity}
\end{align}
is roughly constant.

At each step of the adaptive algorithm we solve the virtual element formulation \eqref{eqn: vem problem}, compute the contribution of each element to the a posteriori error bound in Theorem~\ref{thm: upper bound}, and refine the elements in $\mesh$ with the largest contribution for the next iteration. In order to mark the elements with the largest error contribution, we employ a D\"orfler marking strategy \cite{doerfler}; i.e., we construct the smallest subset of elements $\mesh^M \subset \mesh$  such that
\begin{align*}
    \left( \sum_{\element \in \mesh^M} \eta^2_{\element} + \Theta^2_{\element} + \mathcal{S}^2_{\element} + \Psi^2_{\element} \right)^{\half}
    \geq \theta \left( \sum_{\element \in \mesh} \eta^2_{\element} + \Theta^2_{\element} + \mathcal{S}^2_{\element} + \Psi^2_{\element} \right)^{\half},
\end{align*}
for a steering parameter $\theta \in (0,1)$, by iteratively adding the element with the largest error contribution. For all our numerical experiments we use $\theta=0.4$. In order to refine the elements we note that in our numerical experiments the elements are always convex and, hence, we can refine by connecting the midpoint of each planar edge to the element barycentre; cf. the mVEM package \cite{yu2021implementation} for details on a MATLAB implementation. We note that for non-convex elements obeying Assumption~\ref{ass: mesh regularity} we can use a point to which the element is star-shaped instead.
This refinement strategy introduces hanging nodes; however, these are handled automatically within the VEM framework. We note that there are other refinement strategies available, see e.g. \cite{antonietti2022machine,berrone2021refinement,van2022mesh}.

\pgfplotscreateplotcyclelist{mylist}{
  {black,mark=o},
  {black,mark=triangle},
  {black,mark=square},
  {black,mark=o,densely dashed,every mark/.append style={solid}},
  {black,mark=triangle,densely dashed,every mark/.append style={solid}},
  {black,mark=square,densely dashed,every mark/.append style={solid}},
  {loosely dotted,mark=+},
  {brown!60!black,mark options={fill=brown!40},mark=otimes*}
}

\subsection{Problem 1: smooth solution}\label{subsec: numerics prob4}
\begin{figure}[p]
  \centering
  \subfloat[Initial grid]{
      \includegraphics[trim={3.5cm 0.5cm 3cm 0.5cm},clip, width=0.265\textwidth]{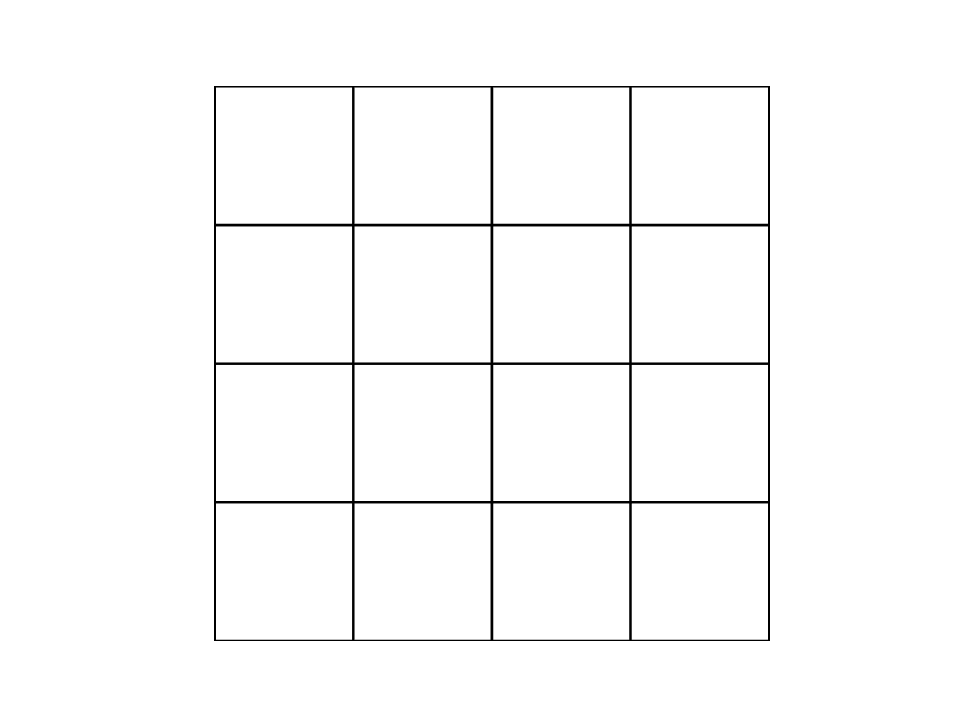}
      \label{fig: prob1 quads step 0}
    }
  \qquad
  \subfloat[After 11 refinements]{
    \includegraphics[trim={3.5cm 0.5cm 3cm 0.5cm},clip, width=0.265\textwidth]{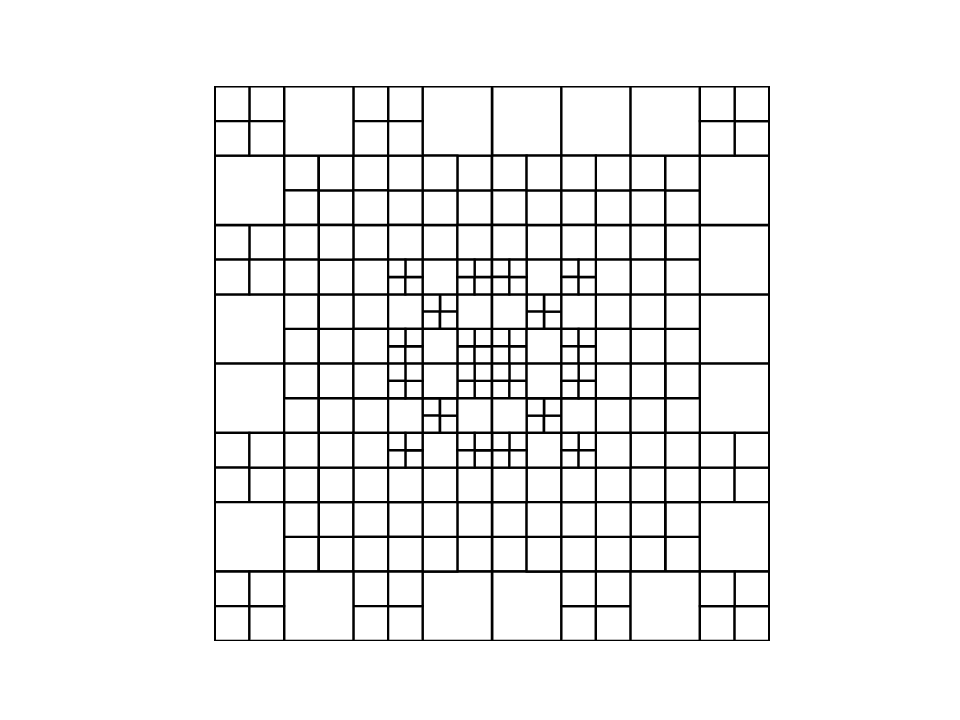}
    \label{fig: prob1 quads step 11}
  }
  \qquad
  \subfloat[After 22 refinements]{
    \includegraphics[trim={3.5cm 0.5cm 3cm 0.5cm},clip, width=0.265\textwidth]{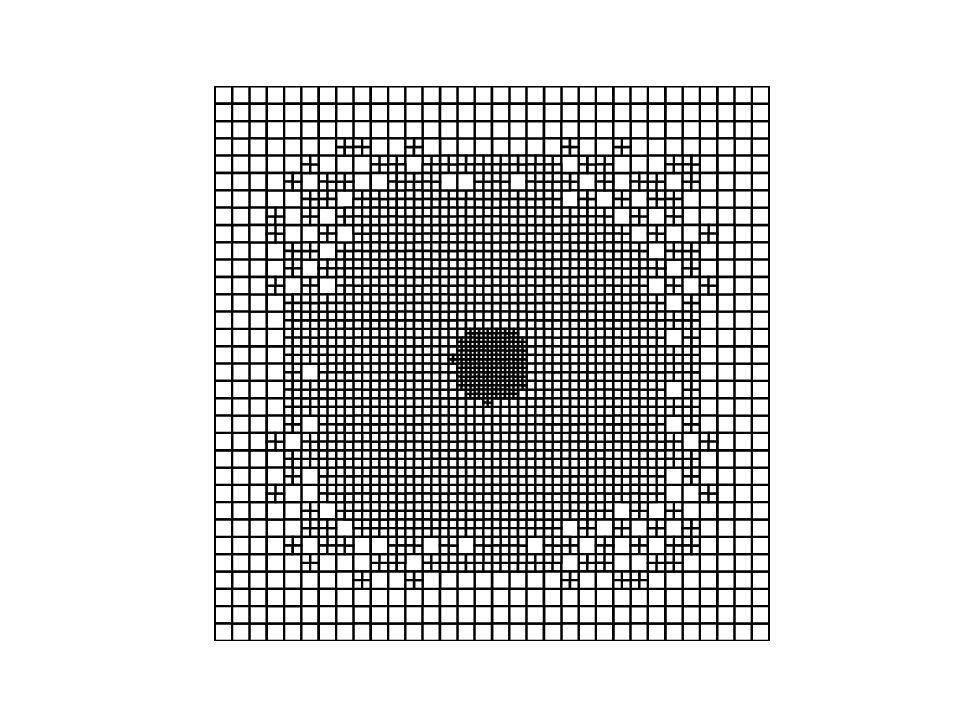}
    \label{fig: prob1 quads step 22}
  }
  \caption{Problem 1: three mesh steps from the adaptive refinement of problem 1 (section~\ref{subsec: numerics prob4}) for the lowest order VEM ($\polOrder=1$) with initial quadrilateral mesh.}
  \label{fig: prob4 quads grids}
\end{figure}

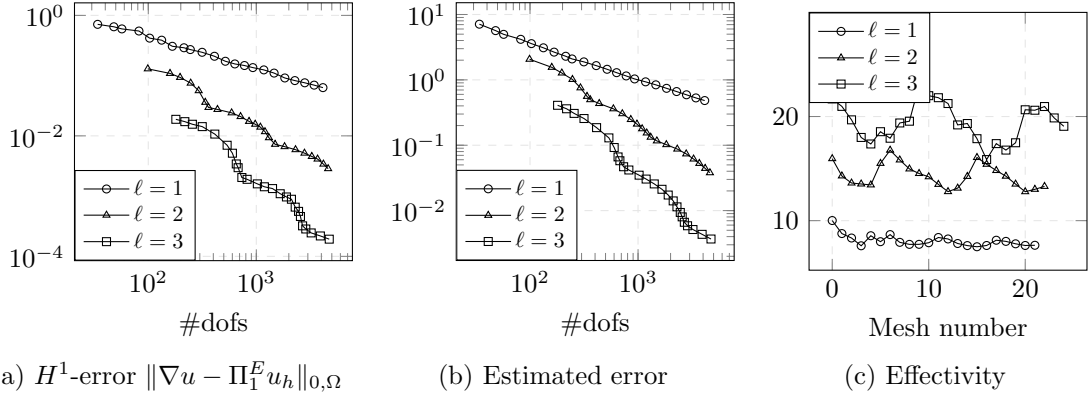
\begin{figure}[p]
    \subfloat[$H^1$-error $\| \nabla u - \gradProj u_h \|_{0,\Omega}$ ]{
      \begin{tikzpicture}
        \begin{loglogaxis}[
            height=5cm,
            width=149pt,
            label style={font=\small},
            tick label style={font=\small},
            grid=major, 
            grid style={dashed,gray!20}, 
            xlabel={\#dofs}, 
            legend style={at=({0,0}), anchor=south west, nodes={scale=0.75, transform shape}},
            mark size=1.5,
            cycle list name=mylist,
            scaled ticks=true
          ]
          \addplot+[] table[x=dofs,y=H1 error,col sep=comma] {revised_results/problem4/problem4_voronoi_p1.csv};
          \addplot+[] table[x=dofs,y=H1 error,col sep=comma] {revised_results/problem4/problem4_voronoi_p2.csv};
          \addplot+[] table[x=dofs,y=H1 error,col sep=comma] {revised_results/problem4/problem4_voronoi_p3.csv};
          \legend{$\polOrder=1$,$\polOrder=2$,$\polOrder=3$}
        \end{loglogaxis}
      \end{tikzpicture}
      \label{fig: voronoi h1 error}}
    \
    \subfloat[Estimated error]{
      \begin{tikzpicture}
        \begin{loglogaxis}[
            height=5cm,
            width=149pt, 
            label style={font=\small},
            tick label style={font=\small},
            grid=major, 
            grid style={dashed,gray!20}, 
            scaled ticks=true,
            scaled y ticks = base 10:1,
            xlabel={\#dofs}, 
            legend style={at=({0,0}), anchor=south west, nodes={scale=0.75, transform shape}},
            mark phase=0,
            mark size=1.5,
            cycle list name=mylist,
          ]
          \addplot+[] table[x=dofs,y=Estimated error,col sep=comma] {revised_results/problem4/problem4_voronoi_p1.csv};
          \addplot+[] table[x=dofs,y=Estimated error,col sep=comma] {revised_results/problem4/problem4_voronoi_p2.csv};
          \addplot+[] table[x=dofs,y=Estimated error,col sep=comma] {revised_results/problem4/problem4_voronoi_p3.csv};
          \legend{$\polOrder=1$,$\polOrder=2$,$\polOrder=3$}
        \end{loglogaxis}
      \end{tikzpicture}
      \label{fig: voronoi estimator error}}
    \
    \subfloat[Effectivity]{
      \begin{tikzpicture}
        \begin{axis}[
            height=5cm,
            width=149pt, 
            label style={font=\small},
            tick label style={font=\small},
            grid=major, 
            grid style={dashed,gray!20}, 
            xlabel={Mesh number}, 
            ymax=30,
            legend style={at=({0,1}), anchor=north west, nodes={scale=0.75, transform shape}},
            mark phase=0,
            mark size=1.5,
            cycle list name=mylist,
            ytick={10, 20},
            yticklabels={$10$,$20$}
          ]
          \addplot+[] table[x=level,y=Effectivity,col sep=comma] {revised_results/problem4/problem4_voronoi_p1.csv};
          \addplot+[] table[x=level,y=Effectivity,col sep=comma] {revised_results/problem4/problem4_voronoi_p2.csv};
          \addplot+[] table[x=level,y=Effectivity,col sep=comma] {revised_results/problem4/problem4_voronoi_p3.csv};
          \legend{$\polOrder=1$,$\polOrder=2$,$\polOrder=3$}
        \end{axis}
      \end{tikzpicture}
      \label{fig: voronoi effectivity}
    }
  \caption{Problem 1: results from solving problem 1 (section~\ref{subsec: numerics prob4}) on the Voronoi grid with adaptive refinement showing (\subref{fig: voronoi h1 error})~convergence history in the $H^1(\Omega)$ seminorm, (\subref{fig: voronoi estimator error})~estimated error, and (\subref{fig: voronoi effectivity}) effectivity \eqref{eqn: effectivity} of the estimator.}
  \label{fig: prob4 voronoi graphs}
\end{figure}

\begin{figure}[p]
  \centering
  \subfloat[Initial grid]{
      \includegraphics[trim={3.5cm 0.5cm 3cm 0.5cm},clip, width=0.25\textwidth]{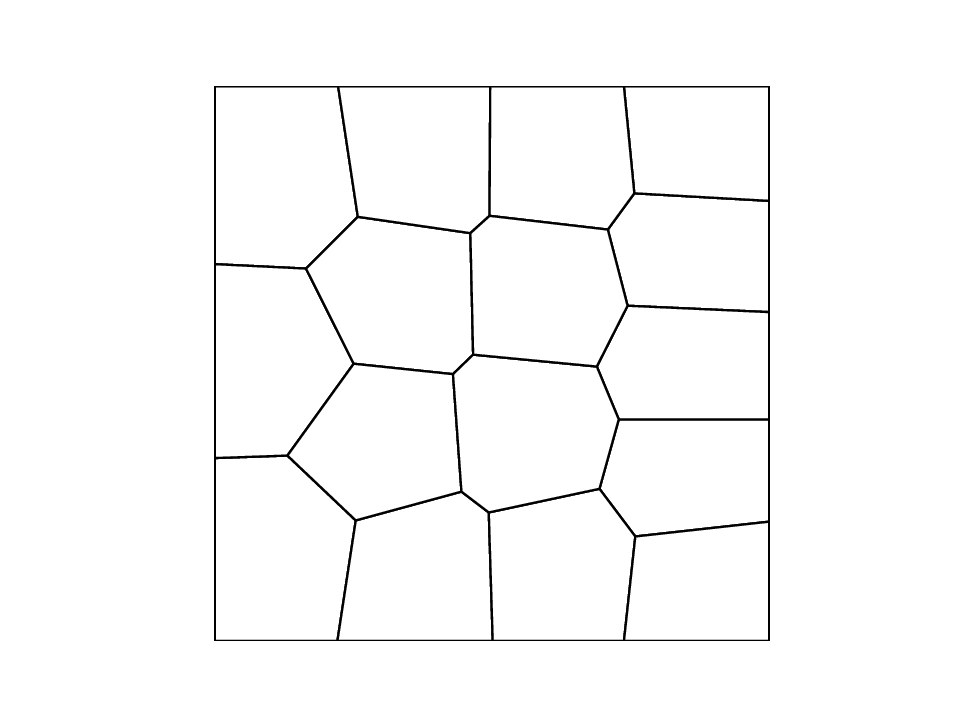}
      \label{fig: prob1 voronoi step 0}
    }
  \qquad
  \subfloat[After 11 refinements]{
    \includegraphics[trim={3.5cm 0.5cm 3cm 0.5cm},clip, width=0.25\textwidth]{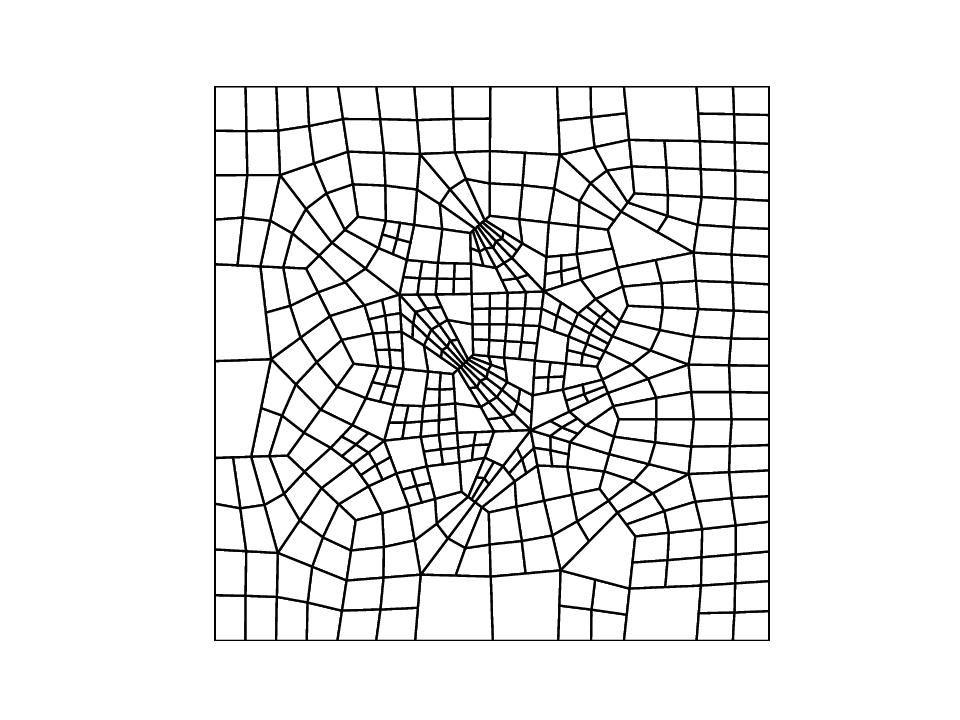}
    \label{fig: prob1 voronoi step 11}
  }
  \qquad
  \subfloat[After 22 refinements]{
    \includegraphics[trim={3.5cm 0.5cm 3cm 0.5cm},clip, width=0.25\textwidth]{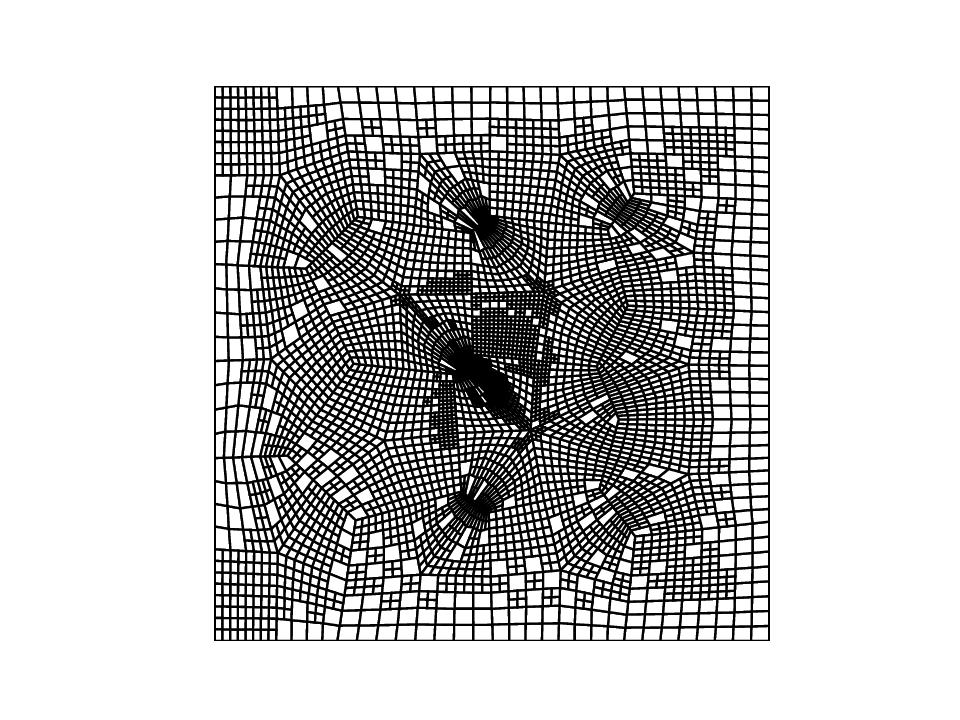}
    \label{fig: prob1 voronoi step 22}
  }
  \caption{Problem 1: three mesh steps from the adaptive refinement of problem 1 (section~\ref{subsec: numerics prob4}) for the lowest order VEM ($\polOrder=1$) with initial Voronoi mesh.}
  \label{fig: prob4 voronoi grids}
\end{figure}

In this first example  we let $\Omega = (0,1)^2$ and we define the nonlinear coefficient $\mu$ as follows
\begin{align*}
    \mu(\boldsymbol{x},|\nabla u|) = 2 + \frac{1}{1+|\nabla u|^2}.
\end{align*}
Furthermore, we take the right hand side $f$ so that the exact solution is given by
\begin{align*}
    u(x,y) = \sin(\pi x)\sin(\pi y).
\end{align*}
We run the adaptive algorithm for the fixed order of approximation $\polOrder=1,2,3$ on both a structured ($4\times4$) quadrilateral grid and a smoothed Voronoi grid (of $16$ elements); cf. Figure~\ref{fig: prob1 quads step 0} and Figure~\ref{fig: prob1 voronoi step 0} for the initial quadrilateral and Voronoi grids, respectively.
In Figures~\ref{fig: quads h1 error}--\ref{fig: quad estimator error} and Figures~\ref{fig: voronoi h1 error}--\ref{fig: voronoi estimator error}, we compare the actual error of the gradient projection of the VEM solution and its a posteriori error bound from Corollary~\ref{cor: upper bound proj} to the number of degrees of freedom for the quadrilateral and Voronoi meshes, respectively, for the sequence of meshes generated by the adaptive mesh refinement algorithm.
We observe that for both initial meshes the actual error and the error bound converges at a similar rate, with the error bound overestimating the true error by a roughly consistent factor; which is confirmed by Figures~\ref{fig: quad effectivity} and~\ref{fig: voronoi effectivity} which show that the effectivity index \eqref{eqn: effectivity} for each mesh is roughly constant, although dependent on the approximation order. We remark that for virtual element methods the effectivity does appear slightly more oscillatory with a larger variance than usual; cf., for example, the effectivity indices for a discontinuous Galerkin finite element method \cite{houston2008posteriori}.

Figures~\ref{fig: prob1 quads step 11}--\ref{fig: prob1 quads step 22} and~\ref{fig: prob1 voronoi step 11}--\ref{fig: prob1 voronoi step 22} display the mesh after 11 and 22 refinements for both the initial quadrilateral and Voronoi meshes, respectively, for $\polOrder=1$. We note that the refinement is roughly uniform, with a slightly more concentrated refinement around the centre, as would be expected for the smooth analytical function considered.

\begin{figure}[p]
    \subfloat[$H^1$-error $\| \nabla u - \gradProj u_h \|_{0,\Omega}$ ]{
      \begin{tikzpicture}
        \begin{loglogaxis}[
            height=5cm,
            width=149pt,
            label style={font=\small},
            tick label style={font=\small},
            grid=major, 
            grid style={dashed,gray!20}, 
            xlabel={\#dofs}, 
            legend style={at=({0,0}), anchor=south west, nodes={scale=0.75, transform shape}},
            mark size=1.5,
            cycle list name=mylist,
            scaled ticks=true
          ]
          \addplot+[] table[x=dofs,y=H1 error,col sep=comma] {revised_results/problem4/problem4_quads_p1.csv};
          \addplot+[] table[x=dofs,y=H1 error,col sep=comma] {revised_results/problem4/problem4_quads_p2.csv};
          \addplot+[] table[x=dofs,y=H1 error,col sep=comma] {revised_results/problem4/problem4_quads_p3.csv};
          \legend{$\polOrder=1$,$\polOrder=2$,$\polOrder=3$}
        \end{loglogaxis}
      \end{tikzpicture}
      \label{fig: quads h1 error}}
    \
    \subfloat[Estimated error]{
      \begin{tikzpicture}
        \begin{loglogaxis}[
            height=5cm,
            width=149pt, 
            label style={font=\small},
            tick label style={font=\small},
            grid=major, 
            grid style={dashed,gray!20}, 
            scaled ticks=true,
            scaled y ticks = base 10:1,
            xlabel={\#dofs}, 
            legend style={at=({0,0}), anchor=south west, nodes={scale=0.75, transform shape}},
            mark phase=0,
            mark size=1.5,
            cycle list name=mylist,
          ]
          \addplot+[] table[x=dofs,y=Estimated error,col sep=comma] {revised_results/problem4/problem4_quads_p1.csv};
          \addplot+[] table[x=dofs,y=Estimated error,col sep=comma] {revised_results/problem4/problem4_quads_p2.csv};
          \addplot+[] table[x=dofs,y=Estimated error,col sep=comma] {revised_results/problem4/problem4_quads_p3.csv};
          \legend{$\polOrder=1$,$\polOrder=2$,$\polOrder=3$}
        \end{loglogaxis}
      \end{tikzpicture}
      \label{fig: quad estimator error}}
    \
    \subfloat[Effectivity]{
      \begin{tikzpicture}
        \begin{axis}[
            height=5cm,
            width=149pt, 
            label style={font=\small},
            tick label style={font=\small},
            grid=major, 
            grid style={dashed,gray!20}, 
            xlabel={Mesh number}, 
            ymax=30,
            legend style={at=({0,1}), anchor=north west, nodes={scale=0.75, transform shape}},
            mark phase=0,
            mark size=1.5,
            cycle list name=mylist,
            ytick={10, 20},
            yticklabels={$10$,$20$}
          ]
          \addplot+[] table[x=level,y=Effectivity,col sep=comma] {revised_results/problem4/problem4_quads_p1.csv};
          \addplot+[] table[x=level,y=Effectivity,col sep=comma] {revised_results/problem4/problem4_quads_p2.csv};
          \addplot+[] table[x=level,y=Effectivity,col sep=comma] {revised_results/problem4/problem4_quads_p3.csv};
          \legend{$\polOrder=1$,$\polOrder=2$,$\polOrder=3$}
        \end{axis}
      \end{tikzpicture}
      \label{fig: quad effectivity}
    }
  \caption{Problem 1: results from solving problem 1 (section~\ref{subsec: numerics prob4}) on the quadrilateral grid with adaptive refinement showing (\subref{fig: quads h1 error})~convergence history in the $H^1(\Omega)$ seminorm, (\subref{fig: quad estimator error})~estimated error, and (\subref{fig: quad effectivity}) effectivity \eqref{eqn: effectivity} of the estimator.
  }
  \label{fig: prob4 quads graphs}
\end{figure}

\FloatBarrier

\subsection{Problem 2: singular solution}\label{subsec: numerics prob2}
In this example, which is also considered in e.g. \cite{wihler2003exponential}, we take $\Omega$ to be the L-shaped domain $\Omega = (-1,1)^2 \, \backslash \, [0,1) \times (-1,0]$ and choose the nonlinear coefficient to be
\begin{align}\label{eqn: prob2 mu}
    \mu(\boldsymbol{x},|\nabla u|) = 1 + e^{- | \nabla u |^2}.
\end{align}
In this case, we choose the forcing $f$ so that the exact solution to \eqref{eqn: pde} is given by
\begin{align*}
    u(r,\theta) = r^{\nicefrac{2}{3}}\sin \left( 2 \theta / 3 \right)
\end{align*}
where $(r,\theta)$ are the usual polar coordinates centred around $(0,0)$. We note that here we additionally impose inhomogeneous Dirichlet boundary conditions.
It is worth noting that in this example, $u$ is analytic in $\bar{\Omega} \backslash \{ \boldsymbol{0} \}$, but $\nabla u$ is singular at the origin.

\begin{figure}[p]
  \centering
  \subfloat[Initial grid]{
      \includegraphics[trim={3.5cm 0.5cm 3cm 0.5cm},clip, width=0.265\textwidth]{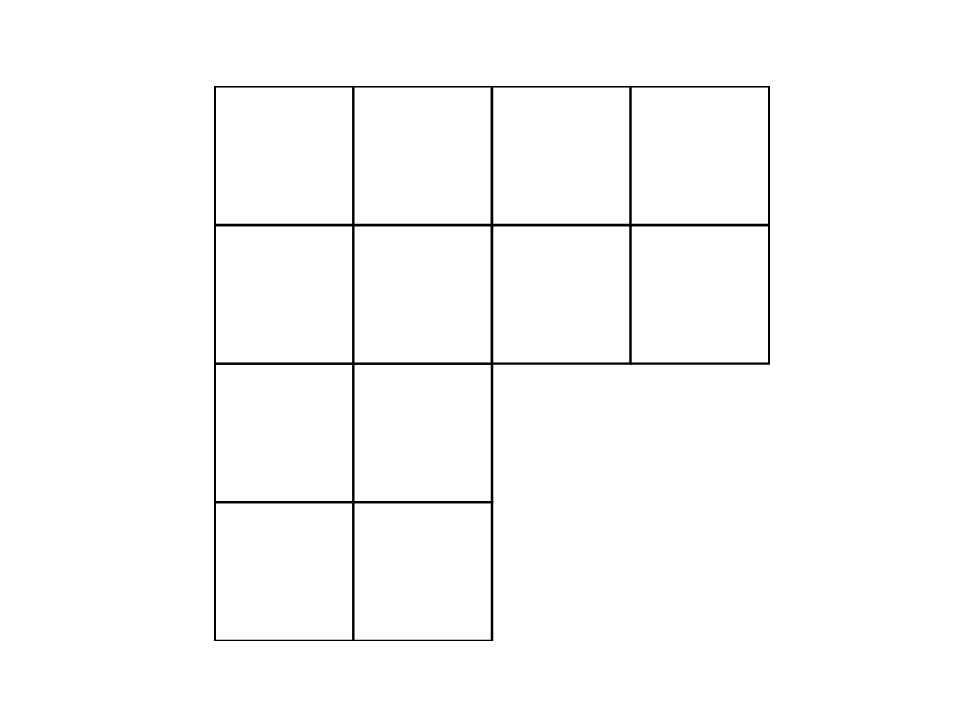}
      \label{fig: prob2 quads step 0}
    }
  \qquad
  \subfloat[After 15 refinements]{
    \includegraphics[trim={3.5cm 0.5cm 3cm 0.5cm},clip, width=0.265\textwidth]{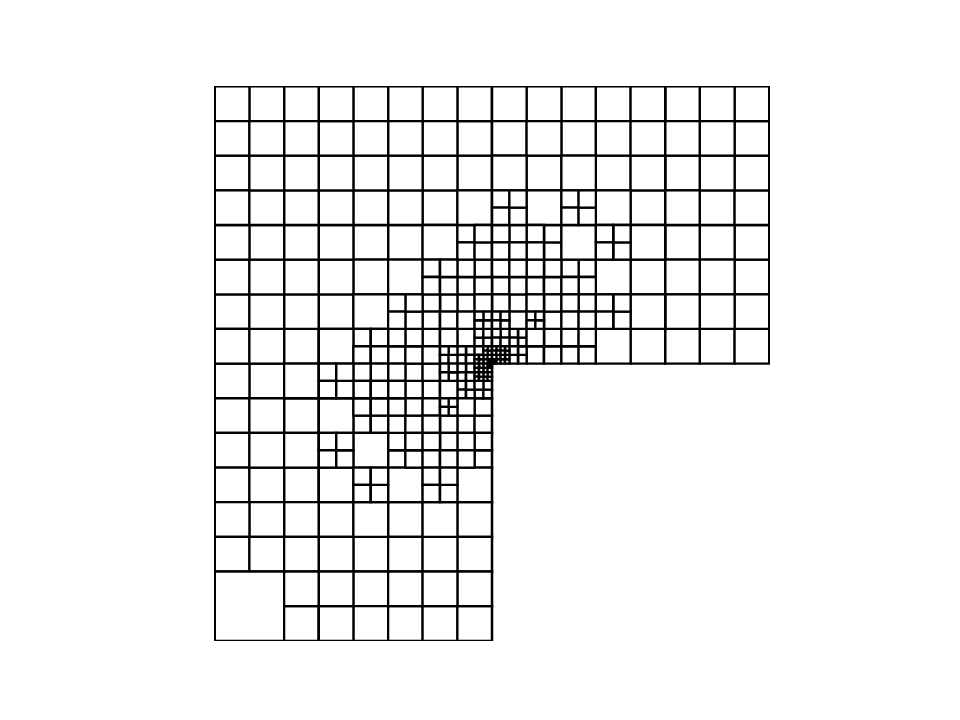}
    \label{fig: prob2 quads step 15}
  }
  \qquad
  \subfloat[After 23 refinements]{
    \includegraphics[trim={3.5cm 0.5cm 3cm 0.5cm},clip, width=0.265\textwidth]{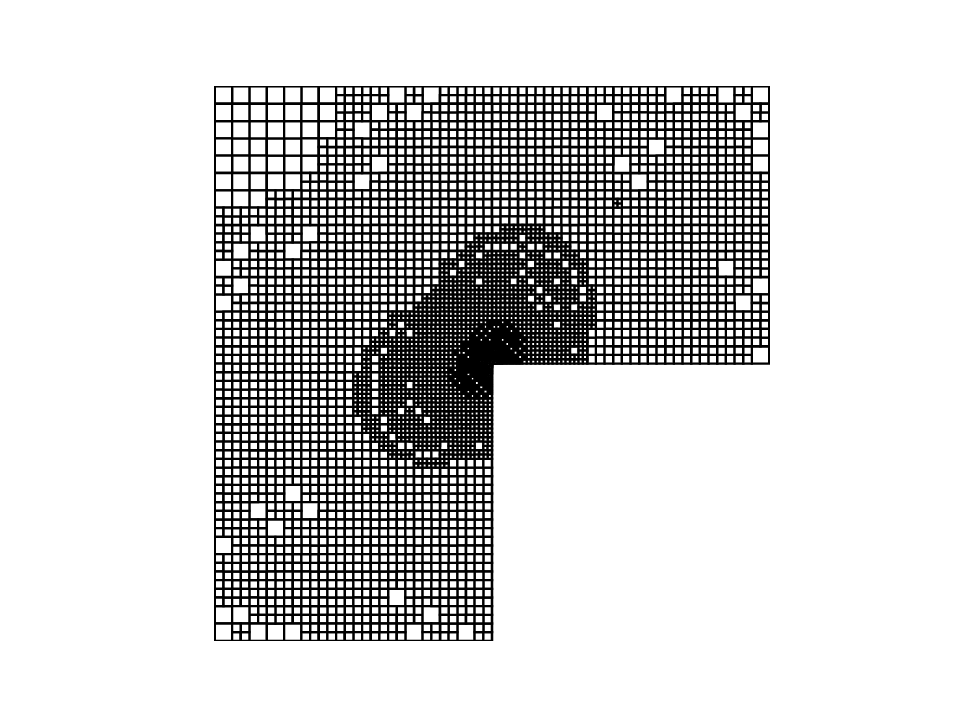}
    \label{fig: prob2 quads step 26}
  }
  \caption{Problem 2: three mesh steps from the adaptive refinement of problem 2 (section~\ref{subsec: numerics prob2}) for the lowest order VEM ($\polOrder=1$) with initial quadrilateral mesh.}
  \label{fig: prob2 quads grids}
\end{figure}

\begin{figure}[p]
  \centering
  \subfloat[Initial grid]{
      \includegraphics[trim={3.5cm 0.5cm 3cm 0.5cm},clip, width=0.265\textwidth]{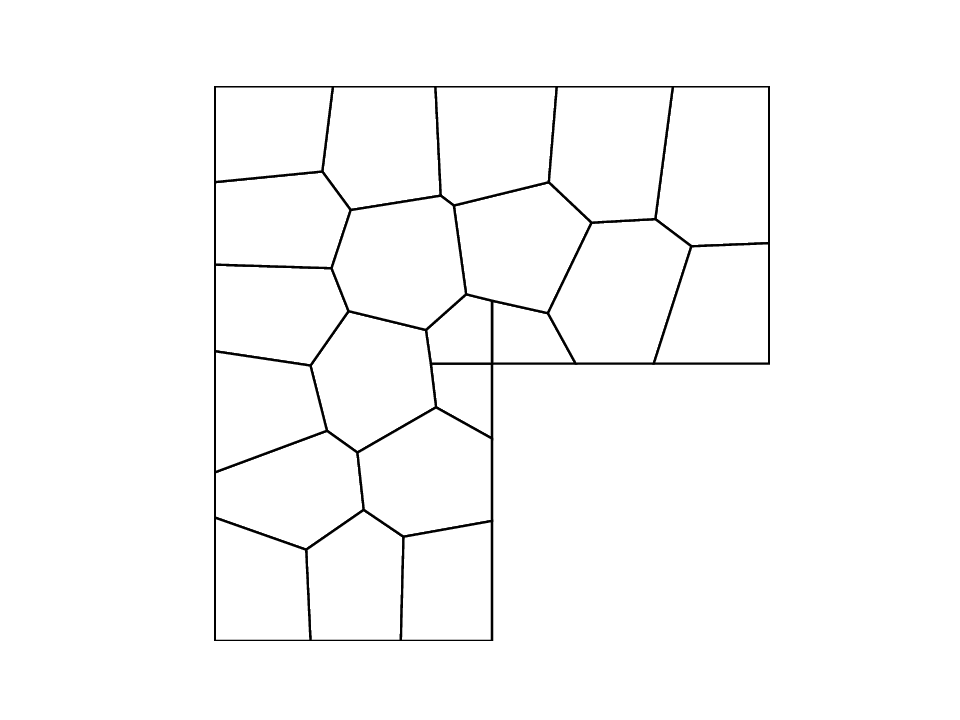}
      \label{fig: prob2 voronoi step 0}
    }
  \qquad
  \subfloat[After 15 refinements]{
    \includegraphics[trim={3.5cm 0.5cm 3cm 0.5cm},clip, width=0.265\textwidth]{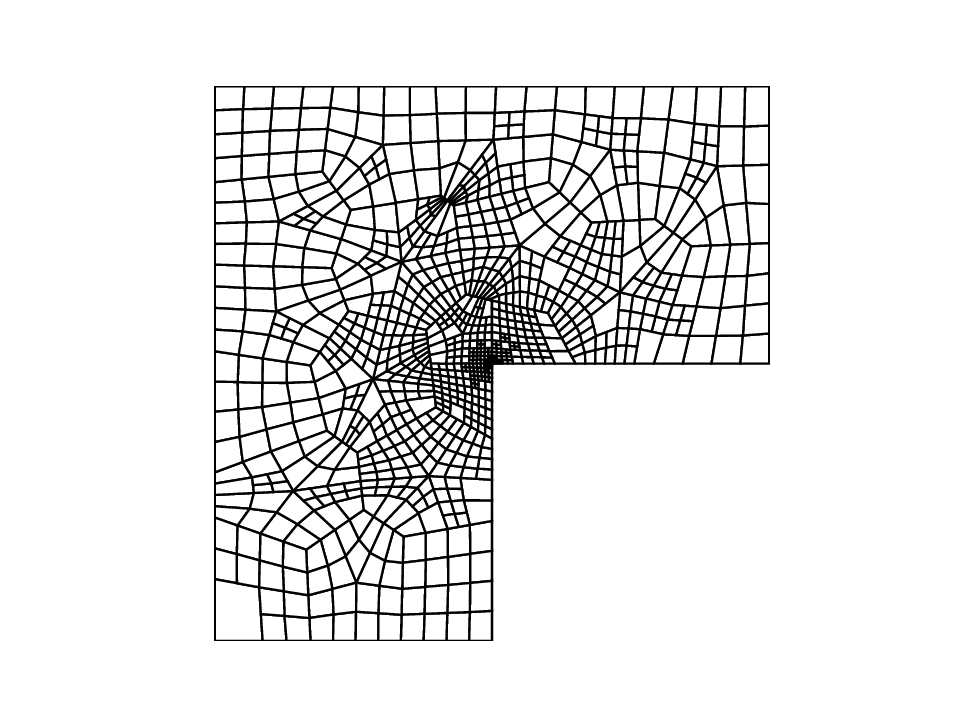}
    \label{fig: prob2 voronoi step 15}
  }
  \qquad
  \subfloat[After 23 refinements]{
    \includegraphics[trim={3.5cm 0.5cm 3cm 0.5cm},clip, width=0.265\textwidth]{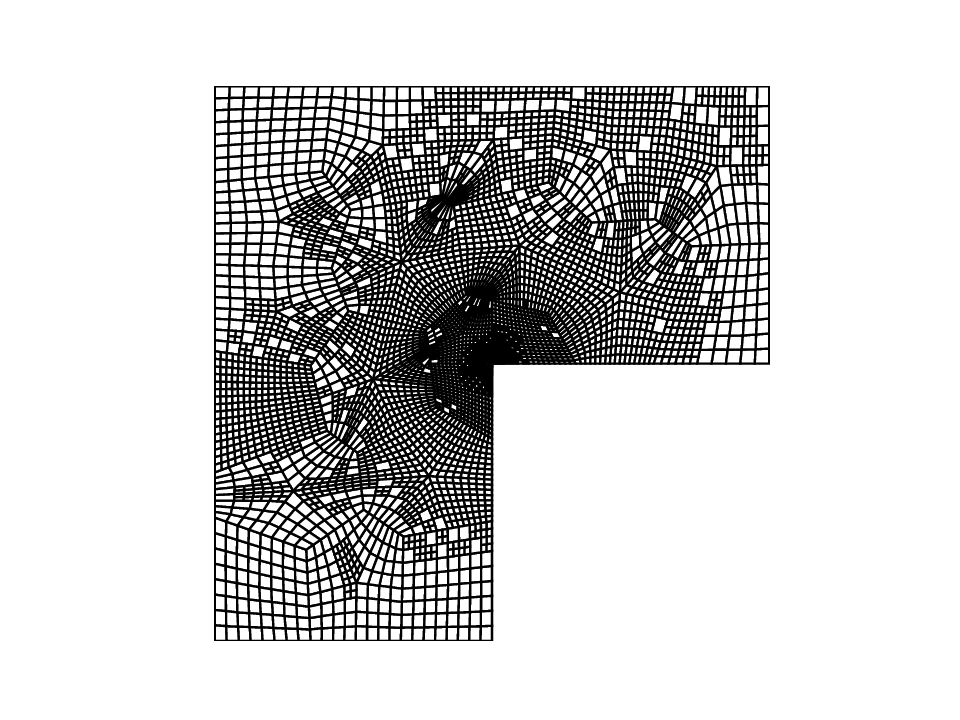}
    \label{fig: prob2 voronoi step 23}
  }
  \caption{Problem 3: three mesh steps from the adaptive refinement of problem 2 (section~\ref{subsec: numerics prob2}) for the lowest order VEM ($\polOrder=1$) with initial Voronoi mesh.}
  \label{fig: prob2 voronoi grids}
\end{figure}

We again run the adaptive algorithm for fixed order of approximation $\polOrder=1,2,3$ on both a structured quadrilateral grid of $12$ elements (Figure~\ref{fig: prob2 quads step 0}) and a smoothed Voronoi grid of $21$ elements (Figure~\ref{fig: prob2 voronoi step 0}). Figures~\ref{fig: prob2 quads h1 error}--\ref{fig: prob2 quads estimator error} and Figures~\ref{fig: prob2 voronoi h1 error}--\ref{fig: prob2 voronoi estimator error} compare the actual error and its a posteriori error bound to the number of degrees of freedom for the quadrilateral and Voronoi meshes, respectively. Again, the actual error and the error bound appear to converge at a similar rate which is confirmed by the effectivity index in Figures~\ref{fig: prob2 quad effectivity} and~\ref{fig: prob2 voronoi effectivity} which appears roughly constant for all meshes.
\begin{figure}[p]
    \subfloat[$H^1$-error $\| \nabla u - \gradProj u_h \|_{0,\Omega}$ ]{
      \begin{tikzpicture}
        \begin{loglogaxis}[
            height=5cm,
            width=149pt,
            label style={font=\small},
            tick label style={font=\small},
            grid=major, 
            grid style={dashed,gray!20}, 
            xlabel={\#dofs}, 
            legend style={at=({0,0}), anchor=south west, nodes={scale=0.75, transform shape}},
            mark size=1.5,
            cycle list name=mylist,
            scaled ticks=true
          ]
          \addplot+[] table[x=dofs,y=H1 error,col sep=comma] {new_results/problem2/problem2_quads_p1.csv};
          \addplot+[] table[x=dofs,y=H1 error,col sep=comma] {new_results/problem2/problem2_quads_p2.csv};
          \addplot+[] table[x=dofs,y=H1 error,col sep=comma] {new_results/problem2/problem2_quads_p3.csv};
          \legend{$\polOrder=1$,$\polOrder=2$,$\polOrder=3$}
        \end{loglogaxis}
      \end{tikzpicture}
      \label{fig: prob2 quads h1 error}}
    \
    \subfloat[Estimated error]{
      \begin{tikzpicture}
        \begin{loglogaxis}[
            height=5cm,
            width=149pt, 
            label style={font=\small},
            tick label style={font=\small},
            grid=major, 
            grid style={dashed,gray!20}, 
            scaled ticks=true,
            scaled y ticks = base 10:1,
            xlabel={\#dofs}, 
            legend style={at=({0,0}), anchor=south west, nodes={scale=0.75, transform shape}},
            mark phase=0,
            mark size=1.5,
            cycle list name=mylist,
          ]
          \addplot+[] table[x=dofs,y=Estimated error,col sep=comma] {new_results/problem2/problem2_quads_p1.csv};
          \addplot+[] table[x=dofs,y=Estimated error,col sep=comma] {new_results/problem2/problem2_quads_p2.csv};
          \addplot+[] table[x=dofs,y=Estimated error,col sep=comma] {new_results/problem2/problem2_quads_p3.csv};
          \legend{$\polOrder=1$,$\polOrder=2$,$\polOrder=3$}
        \end{loglogaxis}
      \end{tikzpicture}
      \label{fig: prob2 quads estimator error}}
    \
    \subfloat[Effectivity]{
      \begin{tikzpicture}
        \begin{axis}[
            height=5cm,
            width=149pt, 
            label style={font=\small},
            tick label style={font=\small},
            grid=major, 
            grid style={dashed,gray!20}, 
            xlabel={Mesh number}, 
            ymax=15,
            legend style={at=({1,1}), anchor=north east, nodes={scale=0.75, transform shape}},
            mark phase=0,
            mark size=1.5,
            cycle list name=mylist,
          ]
          \addplot+[] table[x=level,y=Effectivity,col sep=comma] {new_results/problem2/problem2_quads_p1.csv};
          \addplot+[] table[x=level,y=Effectivity,col sep=comma] {new_results/problem2/problem2_quads_p2.csv};
          \addplot+[] table[x=level,y=Effectivity,col sep=comma] {new_results/problem2/problem2_quads_p3.csv};
          \legend{$\polOrder=1$,$\polOrder=2$,$\polOrder=3$}
        \end{axis}
      \end{tikzpicture}
      \label{fig: prob2 quad effectivity}
    }
  \caption{Problem 2: results from solving problem 2 (section~\ref{subsec: numerics prob2}) on the quadrilateral grid with adaptive refinement showing (\subref{fig: prob2 quads h1 error})~convergence history in the $H^1(\Omega)$ seminorm, (\subref{fig: prob2 quads estimator error})~estimated error, and (\subref{fig: prob2 quad effectivity}) effectivity \eqref{eqn: effectivity} of the estimator.}
  \label{fig: prob2 quads graphs}
\end{figure}

\begin{figure}[p]
    \subfloat[$H^1$-error $\| \nabla u - \gradProj u_h \|_{0,\Omega}$ ]{
      \begin{tikzpicture}
        \begin{loglogaxis}[
            height=5cm,
            width=149pt,
            label style={font=\small},
            tick label style={font=\small},
            grid=major, 
            grid style={dashed,gray!20}, 
            xlabel={\#dofs}, 
            legend style={at=({0,0}), anchor=south west, nodes={scale=0.75, transform shape}},
            mark size=1.5,
            cycle list name=mylist,
            scaled ticks=true
          ]
          \addplot+[] table[x=dofs,y=H1 error,col sep=comma] {new_results/problem2/problem2_voronoi_p1.csv};
          \addplot+[] table[x=dofs,y=H1 error,col sep=comma] {new_results/problem2/problem2_voronoi_p2.csv};
          \addplot+[] table[x=dofs,y=H1 error,col sep=comma] {new_results/problem2/problem2_voronoi_p3.csv};
          \legend{$\polOrder=1$,$\polOrder=2$,$\polOrder=3$}
        \end{loglogaxis}
      \end{tikzpicture}
      \label{fig: prob2 voronoi h1 error}}
    \
    \subfloat[Estimated error]{
      \begin{tikzpicture}
        \begin{loglogaxis}[
            height=5cm,
            width=149pt, 
            label style={font=\small},
            tick label style={font=\small},
            grid=major, 
            grid style={dashed,gray!20}, 
            scaled ticks=true,
            scaled y ticks = base 10:1,
            xlabel={\#dofs}, 
            legend style={at=({0,0}), anchor=south west, nodes={scale=0.75, transform shape}},
            mark phase=0,
            mark size=1.5,
            cycle list name=mylist,
          ]
          \addplot+[] table[x=dofs,y=Estimated error,col sep=comma] {new_results/problem2/problem2_voronoi_p1.csv};
          \addplot+[] table[x=dofs,y=Estimated error,col sep=comma] {new_results/problem2/problem2_voronoi_p2.csv};
          \addplot+[] table[x=dofs,y=Estimated error,col sep=comma] {new_results/problem2/problem2_voronoi_p3.csv};
          \legend{$\polOrder=1$,$\polOrder=2$,$\polOrder=3$}
        \end{loglogaxis}
      \end{tikzpicture}
      \label{fig: prob2 voronoi estimator error}}
    \
    \subfloat[Effectivity]{
      \begin{tikzpicture}
        \begin{axis}[
            height=5cm,
            width=149pt, 
            label style={font=\small},
            tick label style={font=\small},
            grid=major, 
            grid style={dashed,gray!20}, 
            xlabel={Mesh number}, 
            ymax=15,
            legend style={at=({1,1}), anchor=north east, nodes={scale=0.75, transform shape}},
            mark phase=0,
            mark size=1.5,
            cycle list name=mylist,
          ]
          \addplot+[] table[x=level,y=Effectivity,col sep=comma] {new_results/problem2/problem2_voronoi_p1.csv};
          \addplot+[] table[x=level,y=Effectivity,col sep=comma] {new_results/problem2/problem2_voronoi_p2.csv};
          \addplot+[] table[x=level,y=Effectivity,col sep=comma] {new_results/problem2/problem2_voronoi_p3.csv};
          \legend{$\polOrder=1$,$\polOrder=2$,$\polOrder=3$}
        \end{axis}
      \end{tikzpicture}
      \label{fig: prob2 voronoi effectivity}
    }
  \caption{Problem 2: results from solving problem 2 (section~\ref{subsec: numerics prob2}) on the Voronoi grid with adaptive refinement showing (\subref{fig: prob2 voronoi h1 error})~convergence history in the $H^1(\Omega)$ seminorm, (\subref{fig: prob2 voronoi estimator error})~estimated error, and (\subref{fig: prob2 voronoi effectivity}) effectivity \eqref{eqn: effectivity} of the estimator.}
  \label{fig: prob2 voronoi graphs}
\end{figure}

For $\ell=1$ an intermediate and final mesh from the adaptive algorithm are shown in Figures~\ref{fig: prob2 quads step 15}--\ref{fig: prob2 quads step 26} and~\ref{fig: prob2 voronoi step 15}--\ref{fig: prob2 voronoi step 23} for the quadrilateral and Voronoi meshes, respectively. As can be seen the refinement is focused around the singularity at the origin, as would be expected for this singular problem.

\subsection{Problem 3: singular solution with a sharp Gaussian}\label{subsec: numerics prob3}
For the final example, we again consider the L-shaped domain $\Omega = (-1,1)^2 \, \backslash \, [0,1) \times (-1,0]$ from the previous problem, cf. \ref{subsec: numerics prob2}, with the same nonlinearity \eqref{eqn: prob2 mu}, again with a singularity at the origin, but also with an additional sharp Gaussian.
This is similar to the problem considered in \cite{cangiani2017posteriori} for a linear problem. Here, we set the forcing function $f$ and inhomogeneous Dirichlet boundary conditions so that the exact solution is given by
\begin{align*}
    u(x,y) = r^{\nicefrac{2}{3}}\sin \left( 2 \theta / 3 \right) + e^{-(1000(x-0.5)^2 +1000(y-0.5)^2 )}.
\end{align*}
with $(r,\theta)$ denoting the usual polar coordinates, and observe the sharp Gaussian at the point $(0.5,0.5)$.
\begin{figure}[p]
    \subfloat[$H^1$-error $\| \nabla u - \gradProj u_h \|_{0,\Omega}$ ]{
      \begin{tikzpicture}
        \begin{loglogaxis}[
            height=5cm,
            width=149pt,
            label style={font=\small},
            tick label style={font=\small},
            grid=major, 
            grid style={dashed,gray!20}, 
            xlabel={\#dofs}, 
            legend style={at=({0,0}), anchor=south west, nodes={scale=0.75, transform shape}},
            mark size=1.5,
            cycle list name=mylist,
            scaled ticks=true
          ]
          \addplot+[] table[x=dofs,y=H1 error,col sep=comma] {new_results/problem3/problem3_quads_p1.csv};
          \addplot+[] table[x=dofs,y=H1 error,col sep=comma] {new_results/problem3/problem3_quads_p2.csv};
          \addplot+[] table[x=dofs,y=H1 error,col sep=comma] {new_results/problem3/problem3_quads_p3.csv};
          \legend{$\polOrder=1$,$\polOrder=2$,$\polOrder=3$}
        \end{loglogaxis}
      \end{tikzpicture}
      \label{fig: prob3 quads h1 error}}
    \
    \subfloat[Estimated error]{
      \begin{tikzpicture}
        \begin{loglogaxis}[
            height=5cm,
            width=149pt, 
            label style={font=\small},
            tick label style={font=\small},
            grid=major, 
            grid style={dashed,gray!20}, 
            scaled ticks=true,
            scaled y ticks = base 10:1,
            xlabel={\#dofs}, 
            legend style={at=({0,0}), anchor=south west, nodes={scale=0.75, transform shape}},
            mark phase=0,
            mark size=1.5,
            cycle list name=mylist,
          ]
          \addplot+[] table[x=dofs,y=Estimated error,col sep=comma] {new_results/problem3/problem3_quads_p1.csv};
          \addplot+[] table[x=dofs,y=Estimated error,col sep=comma] {new_results/problem3/problem3_quads_p2.csv};
          \addplot+[] table[x=dofs,y=Estimated error,col sep=comma] {new_results/problem3/problem3_quads_p3.csv};
          \legend{$\polOrder=1$,$\polOrder=2$,$\polOrder=3$}
        \end{loglogaxis}
      \end{tikzpicture}
      \label{fig: prob3 quad estimator error}}
    \
    \subfloat[Effectivity]{
      \begin{tikzpicture}
        \begin{axis}[
            height=5cm,
            width=149pt, 
            label style={font=\small},
            tick label style={font=\small},
            grid=major, 
            grid style={dashed,gray!20}, 
            xlabel={Mesh number}, 
            ymax=30,
            legend style={at=({1,1}), anchor=north east, nodes={scale=0.75, transform shape}},
            mark phase=0,
            mark size=1.5,
            cycle list name=mylist,
          ]
          \addplot+[] table[x=level,y=Effectivity,col sep=comma] {new_results/problem3/problem3_quads_p1.csv};
          \addplot+[] table[x=level,y=Effectivity,col sep=comma] {new_results/problem3/problem3_quads_p2.csv};
          \addplot+[] table[x=level,y=Effectivity,col sep=comma] {new_results/problem3/problem3_quads_p3.csv};
          \legend{$\polOrder=1$,$\polOrder=2$,$\polOrder=3$}
        \end{axis}
      \end{tikzpicture}
      \label{fig: prob3 quad effectivity}
    }
  \caption{Problem 3: results from solving problem 1 (section~\ref{subsec: numerics prob3}) on the quadrilateral grid with adaptive refinement showing (\subref{fig: prob3 quads h1 error})~convergence history in the $H^1(\Omega)$ seminorm, (\subref{fig: prob3 quad estimator error})~estimated error, and (\subref{fig: prob3 quad effectivity}) effectivity \eqref{eqn: effectivity} of the estimator.}
  \label{fig: prob3 quads graphs}
\end{figure}

\begin{figure}[p]
  \centering
  \subfloat[Initial grid]{
      \includegraphics[trim={3.5cm 0.5cm 3cm 0.5cm},clip, width=0.265\textwidth]{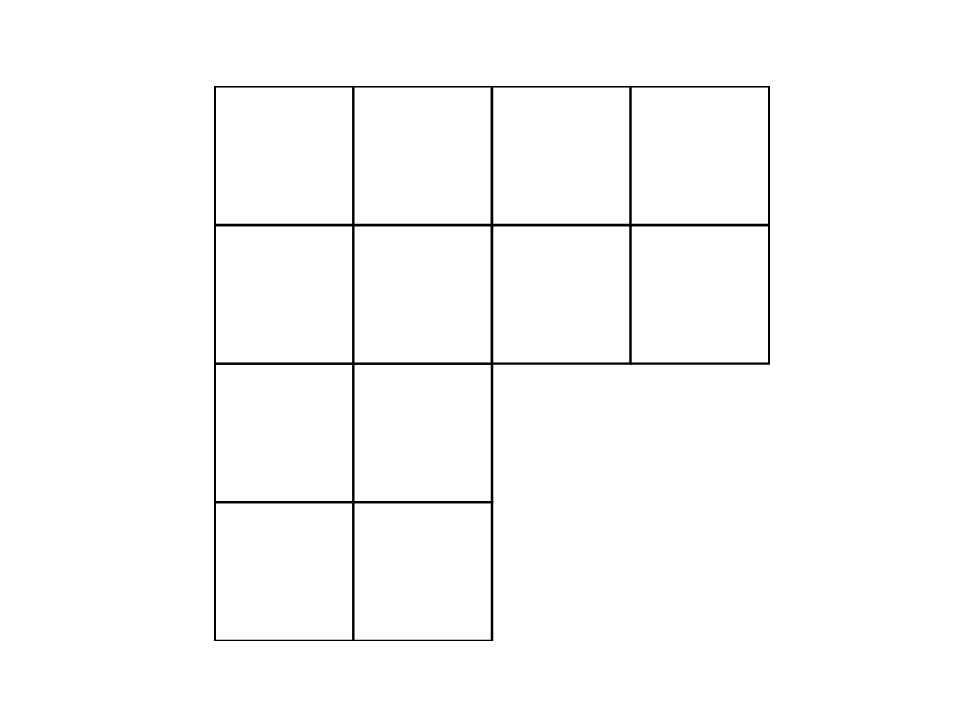}
      \label{fig: prob3 quads step 0}
    }
  \qquad
  \subfloat[After 20 refinements]{
    \includegraphics[trim={3.5cm 0.5cm 3cm 0.5cm},clip, width=0.265\textwidth]{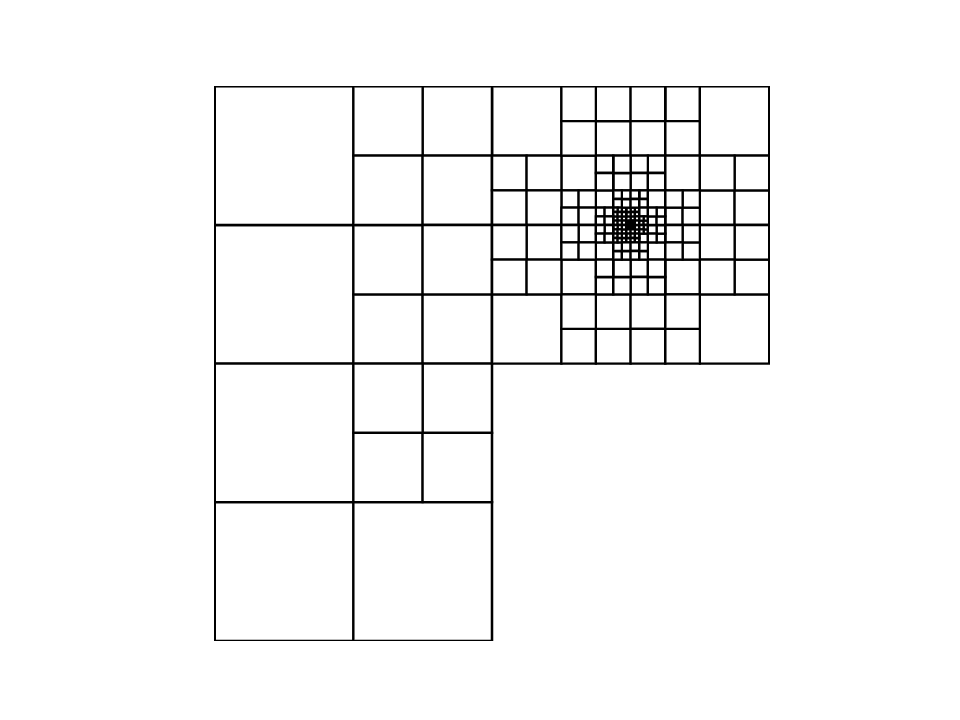}
    \label{fig: prob3 quads step 20}
  }
  \qquad
  \subfloat[After 37 refinements]{
    \includegraphics[trim={3.5cm 0.5cm 3cm 0.5cm},clip, width=0.265\textwidth]{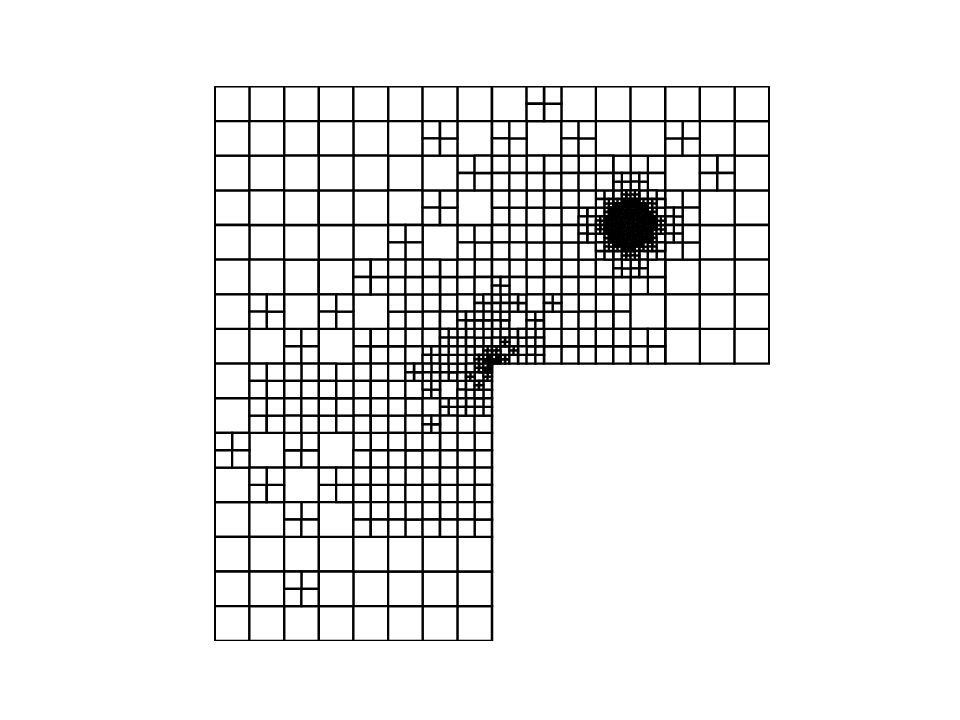}
    \label{fig: prob3 quads step 37}
  }
  \caption{Problem 3: three mesh steps from the adaptive refinement of problem 3 (section~\ref{subsec: numerics prob3}) for the lowest order VEM ($\polOrder=1$) with initial quadrilateral mesh.}
  \label{fig: prob3 quads grids}
\end{figure}

\begin{figure}[p]
    \subfloat[$H^1$-error $\| \nabla u - \gradProj u_h \|_{0,\Omega}$ ]{
      \begin{tikzpicture}
        \begin{loglogaxis}[
            height=5cm,
            width=149pt,
            label style={font=\small},
            tick label style={font=\small},
            grid=major, 
            grid style={dashed,gray!20}, 
            xlabel={\#dofs}, 
            legend style={at=({0,0}), anchor=south west, nodes={scale=0.75, transform shape}},
            mark size=1.5,
            cycle list name=mylist,
            scaled ticks=true
          ]
          \addplot+[] table[x=dofs,y=H1 error,col sep=comma] {new_results/problem3/problem3_voronoi_p1.csv};
          \addplot+[] table[x=dofs,y=H1 error,col sep=comma] {new_results/problem3/problem3_voronoi_p2.csv};
          \addplot+[] table[x=dofs,y=H1 error,col sep=comma] {new_results/problem3/problem3_voronoi_p3.csv};
          \legend{$\polOrder=1$,$\polOrder=2$,$\polOrder=3$}
        \end{loglogaxis}
      \end{tikzpicture}
      \label{fig: prob3 voronoi h1 error}}
    \
    \subfloat[Estimated error]{
      \begin{tikzpicture}
        \begin{loglogaxis}[
            height=5cm,
            width=149pt, 
            label style={font=\small},
            tick label style={font=\small},
            grid=major, 
            grid style={dashed,gray!20}, 
            scaled ticks=true,
            scaled y ticks = base 10:1,
            xlabel={\#dofs}, 
            legend style={at=({0,0}), anchor=south west, nodes={scale=0.75, transform shape}},
            mark phase=0,
            mark size=1.5,
            cycle list name=mylist,
          ]
          \addplot+[] table[x=dofs,y=Estimated error,col sep=comma] {new_results/problem3/problem3_voronoi_p1.csv};
          \addplot+[] table[x=dofs,y=Estimated error,col sep=comma] {new_results/problem3/problem3_voronoi_p2.csv};
          \addplot+[] table[x=dofs,y=Estimated error,col sep=comma] {new_results/problem3/problem3_voronoi_p3.csv};
          \legend{$\polOrder=1$,$\polOrder=2$,$\polOrder=3$}
        \end{loglogaxis}
      \end{tikzpicture}
      \label{fig: prob3 voronoi estimator error}}
    \
    \subfloat[Effectivity]{
      \begin{tikzpicture}
        \begin{axis}[
            height=5cm,
            width=149pt, 
            label style={font=\small},
            tick label style={font=\small},
            grid=major, 
            grid style={dashed,gray!20}, 
            xlabel={Mesh number}, 
            ymax=30,
            legend style={at=({1,1}), anchor=north east, nodes={scale=0.75, transform shape}},
            mark phase=0,
            mark size=1.5,
            cycle list name=mylist,
          ]
          \addplot+[] table[x=level,y=Effectivity,col sep=comma] {new_results/problem3/problem3_voronoi_p1.csv};
          \addplot+[] table[x=level,y=Effectivity,col sep=comma] {new_results/problem3/problem3_voronoi_p2.csv};
          \addplot+[] table[x=level,y=Effectivity,col sep=comma] {new_results/problem3/problem3_voronoi_p3.csv};
          \legend{$\polOrder=1$,$\polOrder=2$,$\polOrder=3$}
        \end{axis}
      \end{tikzpicture}
      \label{fig: prob3 voronoi effectivity}
    }
  \caption{Problem 3: results from solving problem 1 (section~\ref{subsec: numerics prob3}) on the Voronoi grid with adaptive refinement showing (\subref{fig: prob3 voronoi h1 error})~convergence history in the $H^1(\Omega)$ seminorm, (\subref{fig: prob3 voronoi estimator error})~estimated error, and (\subref{fig: prob3 voronoi effectivity}) effectivity \eqref{eqn: effectivity} of the estimator.}
  \label{fig: prob3 voronoi graphs}
\end{figure}

\begin{figure}[p]
  \centering
  \subfloat[Initial grid]{
      \includegraphics[trim={3.5cm 0.5cm 3cm 0.5cm},clip, width=0.265\textwidth]{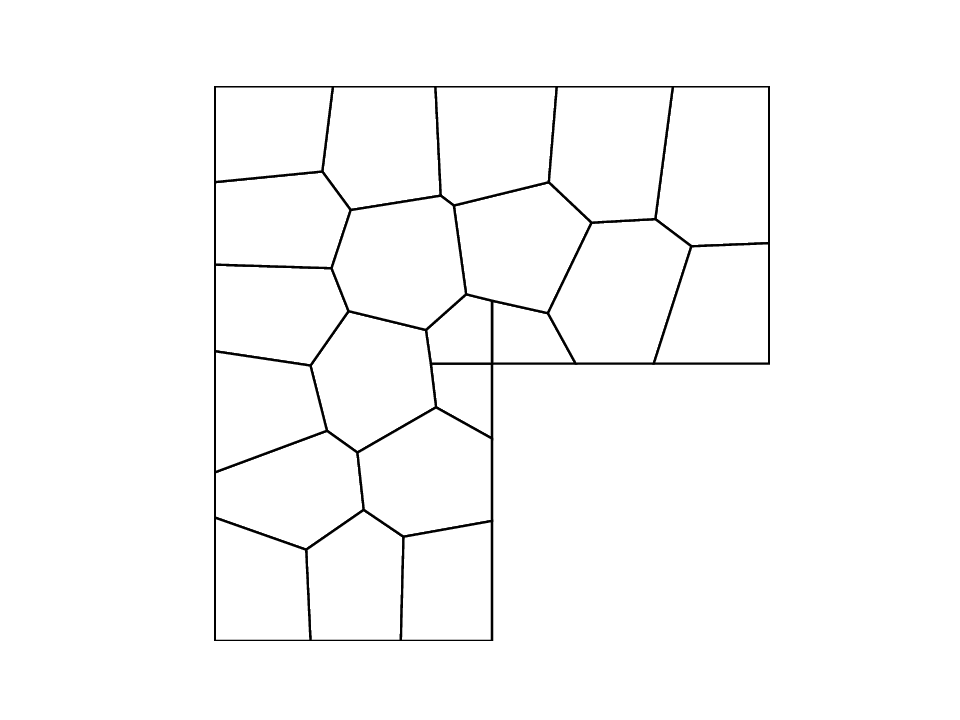}
      \label{fig: prob3 voronoi step 0}
    }
  \qquad
  \subfloat[After 13 refinements]{
    \includegraphics[trim={3.5cm 0.5cm 3cm 0.5cm},clip, width=0.265\textwidth]{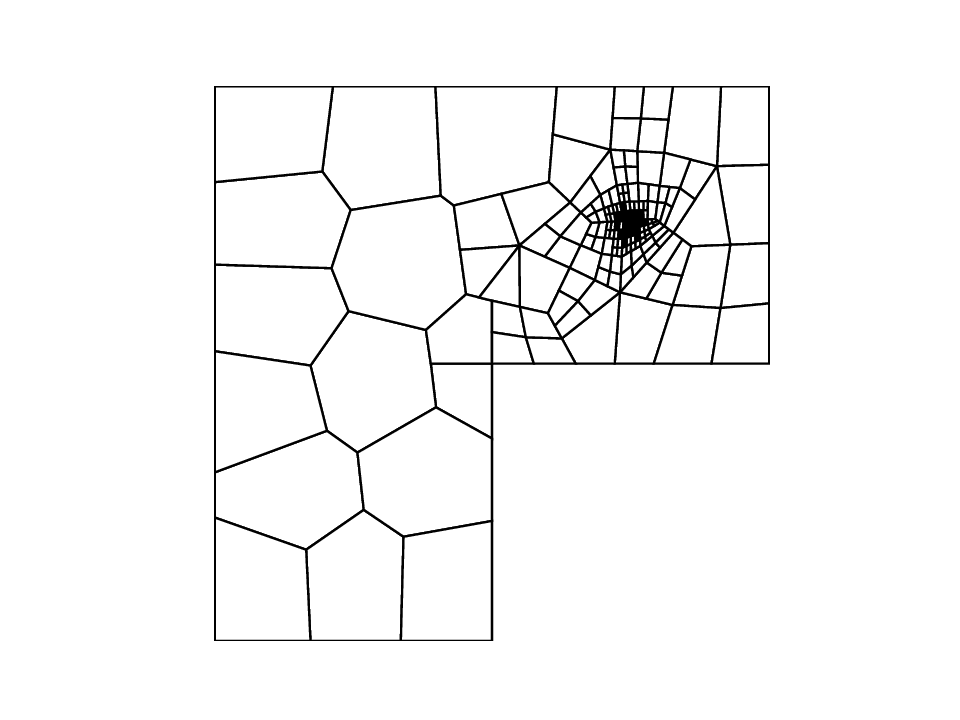}
    \label{fig: prob3 voronoi step 13}
  }
  \qquad
  \subfloat[After 27 refinements]{
    \includegraphics[trim={3.5cm 0.5cm 3cm 0.5cm},clip, width=0.265\textwidth]{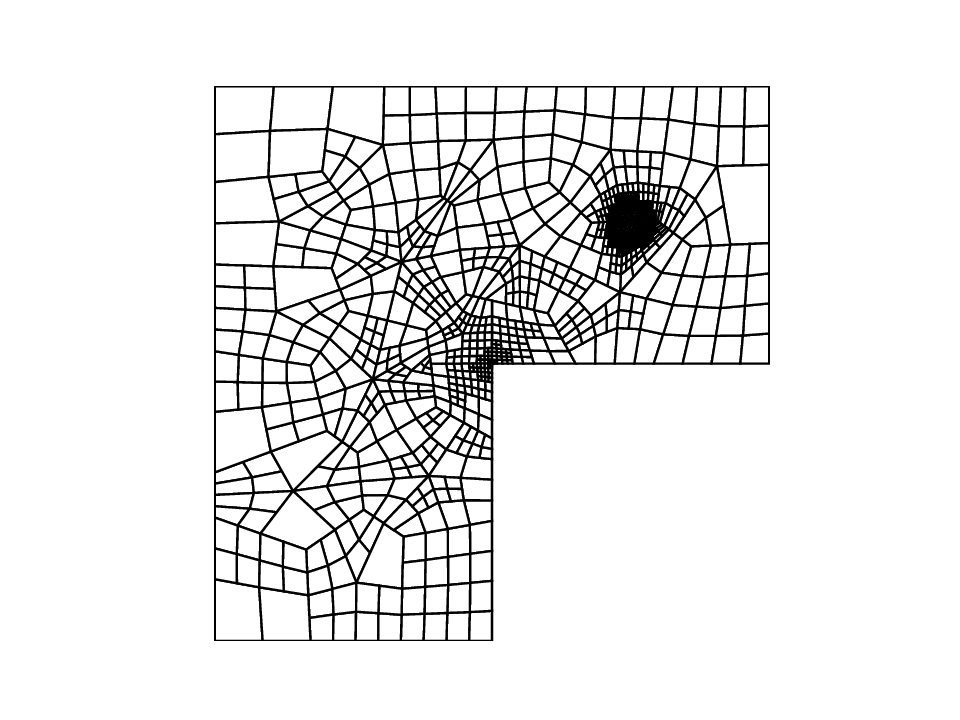}
    \label{fig: prob3 voronoi step 27}
  }
  \caption{Problem 3: three mesh steps from the adaptive refinement of problem 3 (section~\ref{subsec: numerics prob3}) for the lowest order VEM ($\polOrder=1$) with initial Voronoi mesh.}
  \label{fig: prob3 voronoi grids}
\end{figure}

We run the adaptive algorithm for $\polOrder=1,2,3$ on the same initial meshes as the previous problem; cf., Figures~\ref{fig: prob3 quads step 0} and~\ref{fig: prob3 voronoi step 0}. The comparison of the actual error and the error bound is shown in Figures~\ref{fig: prob3 quads h1 error}--\ref{fig: prob3 quad estimator error} and Figures~\ref{fig: prob3 voronoi h1 error}--\ref{fig: prob3 voronoi estimator error}. Here, we see that initial the error bound and true error behave slightly differently, but after some initial pre-asymptotic steps appear to start to converge at a similar rate. This is confirmed by the effectivity indices, cf. Figures~\ref{fig: prob3 quad effectivity} and~\ref{fig: prob3 voronoi effectivity}, which, are roughly constant after some initial steps. We note this initial behaviour likely occurs while the Gaussian is not sufficiently resolved.

An intermediate and final mesh for $\ell=1$ are shown in Figures~\ref{fig: prob3 quads step 20}--\ref{fig: prob3 quads step 37} and~\ref{fig: prob3 voronoi step 13}--\ref{fig: prob3 voronoi step 27}.
Here, it can be seen that the adaptive algorithm first refines around the sharp Gaussian at $(0.5,0.5)$ until it is sufficiently resolved, and then further refinement is focused around this Gaussian and the singularity at the re-entrant corner.

\section{Conclusion}\label{sec: conclusion}
In this paper we have developed a $C^0$-conforming virtual element method of arbitrary approximation order for the discretisation of the second-order quasilinear elliptic PDE in two dimensions \eqref{eqn: pde}.
We have applied the projection approach taken in \cite{10.1093/imanum/drab003,dedner2022framework} and as a result we were able to define the discrete forms directly.
In particular, we discretised the nonlinearity $\mu$ using the gradient projection $\gradProj$ which itself was shown to be the $L^2(\element)$-orthogonal projection of the gradient.
Furthermore, we presented a posteriori error analysis and derived a fully computable residual based error estimator.
Upper and lower bounds for the estimator were shown using techniques from \cite{cangiani2017posteriori} including a standard but important VEM interpolation result, which we detailed in Theorem~\ref{thm: vem approximation}.
Finally, we presented a set of numerical results to study the behaviour of the proposed error estimator when using it to drive an adaptive algorithm.
A variety of tests from the literature were carried out for two different sets of polygonal grids, and we demonstrated that the convergence rate of the a posteriori error bound and true error is roughly similar.

\section*{Funding}
Both authors have been supported by Charles University Research programme no. {PRIMUS/22/SCI/014}.

\bibliographystyle{acm}
\bibliography{quasilinearVEM}

\end{document}